\documentclass{amsart}
 \usepackage[dvips]{color,graphicx}
\usepackage{amsmath,amstext,amssymb,amsbsy,amscd,fancybox}
\usepackage{pst-all}

\newtheorem{defi}{Definition}[section]
\newtheorem{teo}{Theorem}

\newtheorem{prop}{Proposition}
\newtheorem{lema}{Lemma}
\usepackage[active]{srcltx}

\newtheorem{afir}{Claim}[section]

\newtheorem{obs}{Remark}

\newtheorem{conj}{Conjecture}

\newcommand{\pp}{ \vspace{-3mm}\text{}\\ {\bf Proof. }}
\newcommand{\B}{ $\rule{1.2mm}{2mm}$\\}

\begin{document}

\title[ROBUST TRANSITIVITY FOR ENDOMORPHISMS]{ROBUST TRANSITIVITY FOR ENDOMORPHISMS}

\author[C. Lizana]{Cristina $\text{Lizana}^{\dag}$}

\address{$\dag$ Departamento de Matem\'atica\\
Facultad de Ciencias\\
Universidad de Los Andes\\
La Hechicera-M\'{e}rida, 5101\\
Venezuela}
 \email{clizana@ula.ve}
 \thanks{\noindent $\dag$
This work was partially supported by TWAS-CNPq and Universidad de
Los Andes}

\author[E. Pujals]{Enrique  $\text{Pujals}^{\ddag}$}
\address{$\ddag$ IMPA, Estrada Dona Castorina,110, CEP 22460-320, Rio de Janeiro, Brazil}
 \email{enrique@impa.br}
\date{\today}

\maketitle


\begin{abstract}
We address the problem about under what conditions an endomorphism
having a dense orbit, verifies that a sufficiently close perturbed
map also exhibits a dense orbit. In this direction, we give
sufficient conditions, that cover a large class of examples,
for
endomorphisms on the $n-$dimensional torus 
to be robustly transitive: the endomorphism must be volume expanding and any large
connected arc must  contain a point  such that its future orbit
belong to  an expanding region.


\end{abstract}


\section{Introduction}

One goal in dynamics is to look for conditions that guarantee that certain phenomena
is robust under perturbations, that is, under which hypothesis   some main feature of a dynamical system
is shared by all nearby systems. In particular, we are interested in the hypotheses under which
 an endomorphism is robust transitive (see definitions \ref{tran} and \ref{RT}). 

In the diffeomorphism case, there are many examples of robust transitive systems.
The best known is the transitive Anosov diffeomorphism. In the nonhyperbolic context, the first  example
was given by Shub in $\mathbb{T}^4$ in 1971 (see \cite{Shub2}); another example is the
Ma\~{n}\'{e}'s Derived from an Anosov in  $\mathbb{T}^3$ (see \cite{MR});
Bonatti and D\'{\i}az \cite{BD} gave a geometrical construction that produce partially hyperbolic
robust transitive systems and  these constructions were generalized by
Bonatti and Viana providing robust transitive diffeomorphisms with dominated splitting which
are not partially hyperbolic (see \cite{BV}). All those examples are adapted
(and some new ones are extended) to the case of endomorphisms (see section \ref{s5}).

On the other hand, any $C^1-$robust transitive diffeomorphism exhibits a do\-mi\-na\-ted splitting (see \cite{BDP}).
This is no longer true for endomorphisms (see example 1 in section \ref{ss31}).
Therefore, for endomorphisms, conditions that imply robust transitivity cannot hinge on the existence
of splitting.


%

The first question that arises is what necessary condition a robust transitive
 endomorphism has to verify.
Adapting some parts of the proof in \cite{BDP} it is  shown in 
Theorem \ref{teo0}, section \ref{s4}, that for endomorphisms not exhibiting  a dominated
splitting (in a robust way, see definition \ref{nosp}),
volume expanding is a $C^1$ necessary condition.
However, volume expanding is not a sufficient condition that guarantees robust transitivity
for a local diffeomorphism,  as it follows considering  an expanding endomorphism times an irrational rotation (this system
is volume expanding and transitive but not robust transitive, see remark \ref{obs16} for more details). 
Hence, we need an extra condition (that persists by perturbations and does
 not depend on the existence of any type of splitting)
that allow us to conclude robustness.
The extra  hypothesis  that we require can be formulated as follows:
\emph{any arc of diameter large enough   have a point such that its forward
iterates remain in some  expanding region }(see Main Theorem below).



Before introducing the Main Theorem, we recall some definitions  and we  introduce some notation that we  use throughout this work.

An \emph{endomorphism} of a differentiable manifold $M$ is 
a differentiable
function $f:M\rightarrow M$ of class $C^r$ with $r\geq 1.$
Let us denote by $E^r(M)$ ($r\geq 1$) the space of $C^r-$endomorphisms of $M$
endowed with the usual $C^r$ topology. A \emph{local diffeomorphism} is an endomorphism
$f:M\rightarrow M$ such  that given any point $x\in M,$
there exists an open set $V$ in $M$ containing $x$ such that $f$ from
$V$ to $f(V)$ is a diffeomorphism.

\begin{defi}
We say that a map $f\in E^1(M)$ is \emph{volume expanding} if there exists $\sigma>1$
such that $|det(Df)|>\sigma.$
\end{defi}

Observe that volume expanding endomorphisms are
local diffeomorphisms.

If $L:V\rightarrow W$ is a linear isomorphism between
normed vector spaces, we denote by $m\{L\}$ the \emph{minimum norm} of $L,$
 i.e. $m\{L\}=\|L^{-1}\|^{-1}.$

\begin{defi}
We say that a set $\Lambda\subset M$ is a \emph{forward invariant set} for  $f\in E^r(M)$ 
if $f(\Lambda)\subset\Lambda$ and it is \emph{invariant} for $f$ if $f(\Lambda)=\Lambda$.
\end{defi}

\begin{defi}
We say that a map $f\in E^1(M)$  is \emph{expanding} in $U$, a subset of $M,$
if there exists $\lambda>1$
such that $\displaystyle\min_{x\in U}\{m\{D_xf\}\}>\lambda.$ It is said that a compact
invariant set $\Lambda$ is an \emph{expanding set} for an endomorphism $f$ if $f\mid_{\Lambda}$ is
an expanding map.
\end{defi}

\begin{defi}
Let $U$ be an open set in $\mathbb{T}^n$. Denote by $\widetilde{U}$ the lift of $U$
restricted to a fundamental domain of $\mathbb{R}^n$ . 
Define the \emph{ diameter} of $U$ by $$diam(U)= \max\{dist(x, y): x, y\in  \widetilde U\}.$$
Define
the \emph{internal diameter} of $U^c$  by
$$diam_{int}(U^c)\!=\!\displaystyle\min_{k\in\mathbb{Z}^n\setminus \{0\}}dist\!(\widetilde{U}, \widetilde{U}+k),$$
where $dist(A,B):=\displaystyle\inf\{\max_{1\leq i\leq n}|x_i-y_i|: x=(x_1,\ldots,x_n)\!\in\! A, y=(y_1,\ldots,y_n)\!\in\! B\}.$
\end{defi}
Related to the last definition, observe that if $diam(U)<1$ then, translating the frame $\mathbb{Z}^n,$ we can assume that
$\widetilde{U}$ is contained in the interior of $[0,1]^n$ and in particular, $diam_{int}(U^c)>0$.

%

 \begin{defi}\label{tran}
 Let $\Lambda$ be an invariant set for an endomorphism $f:M\rightarrow M.$ It is said
 that $\Lambda$ is 
 \emph{topologically transitive} (or transitive) if there exists a point $x\in \Lambda$
 such that its forward orbit $\{f^k(x)\}_{k\geq 0}$ is dense
 in $\Lambda.$ We say that $f$ is \emph{topologically transitive}
 if $\{f^k(x)\}_{k\geq 0}$ is dense in $M$ for some $x\in M.$
 \end{defi}

The following lemma is a more useful characterization of transitivity.

 \begin{lema}
 Let $f:M\rightarrow M$ be a continuous map of a locally compact
 separable metric space $M$ into itself. The map  $f$ is topologically
 transitive if and only if for any two nonempty open sets
 $U,V\subset M,$ there exists a positive integer $N=N(U,V)$
 such that  $f^N(U)\cap V$ is nonempty.
 \end{lema}

 \pp See for instance \cite[pp.29]{Katok}.\B

Instead of transitivity we may assume the  density of the
pre-orbit of any point. Observe that this implies transitivity,
but the reciprocal assertion does not necessarily hold. In fact,  it is enough
to have a dense subset of the manifold such that every point in this set has dense
pre-orbit to obtain transitivity. On the contrary of diffeomorphisms case,
for endomorphisms,  to have
just one point with dense pre-orbit is not enough to guarantee transitivity.
We leave the details to the reader, it is not hard to construct an example
having some points with dense pre-orbit but non-transitive.

\vspace{5mm}
Let us state the main theorem of the present work.

\vspace{5mm}
\begin{mth}
\emph{Let $f\in E^r(\mathbb{T}^n)$  be a volume expanding map
($n\geq 2, r\geq 1$) such that $\{w\in f^{-k}(x): k\in
\mathbb{N}\}$ is dense for every $x\in \mathbb{T}^n$ and
satisfying the following pro\-per\-ties:
\begin{enumerate}
\item There is an open set $U_0$ in $\mathbb{T}^n$ such that
      $f\!\!\mid_{U_0^c}$ is expanding and  $diam\!(U_0)\!\!<\!\!1$.
\item There exists $0<\delta_0<diam_{int}(U_0^c)$ and there exists an open neighborhood $U_1$
        of $\overline{U}_0$ such that for every  arc $\gamma$ in $U_0^c$ with
      diameter larger than  $\delta_0,$ there is a point $y\in\gamma$ such that
      $f^k(y)\in U_1^c$ for any $k\geq 1.$
\item Moreover, for every $z\in U_1^c,$ there exists $\bar{z}\in U_1^c$
      such that $f(\bar{z})=z.$
\end{enumerate}
Then, for every $g$ $C^r-$close enough to $f,$ $\{w\in g^{-k}(x):
k\in \mathbb{N}\}$ is dense for every $x\in \mathbb{T}^n.$ In
particular, $f$ is $C^r-$robust transitive.}
\end{mth}
\vspace{5mm}

We would like to say a few words about the hypotheses of the Main Theorem. The first
hypothesis states that there exists a set $U_0$ (not necessarily connected) where $f$ fails to be expanding
(if $U_0$ is empty, then $f$ is expanding and the thesis follows from standard arguments for expanding maps),
however, $U_0$ is contained in a ball of radius one and in the complement of it, $f$ is expanding. The second
 hypothesis states that for any large connected arc in the expanding region, there is a point that its forward iterates remains in the expanding region. We assume $n=dim
\mathbb{T}^n$ greater or equal 2, since in dimension 1 if a map is
volume expanding, then it is an expanding map.

A class of systems that verifies the hypotheses of the Main
Theorem is a certain type of maps isotopic to expanding
endomorphisms. More precisely, we call those maps as
``\emph{Derived from Expanding}", the reason to use this name is
inspired on the \emph{Derived from Anosov} (see \cite{MR}) which
are maps isotopic to an Anosov but they are not Anosov. In
particular,
Derived from Expanding maps that satisfies the hypotheses of Main
Theorem  are robustly transitive.
In examples 1 and 2 in section \ref{s5}, 
we show that there exist Derived from Expanding maps
satisfying the hypotheses of Main Theorem.
 We want to point out that in the hypotheses of the Main
Theorem it is not assumed that $f$ is isotopic to an expanding map.


Some questions that arises from the above discussion are: \emph{if
a map satisfies the hypotheses of the Main Theorem, then is 
this map  isotopic to an expanding endomorphism?} \emph{Are robust
transitive endomorphisms
without dominated splitting (in a robust way) isotopic to
expanding endomorphisms?}

We suggest to the reader that before entering into the proof of
Main Theorem, to give a glance to section \ref{ss12} in order to
gain some insight about the proof. We want to highlight that this
theorem as it is enunciated, does not assume the existence of a
tangent bundle splitting (neither it assumes the lack of a
dominated splitting) and it covers examples of robust transitive
endomorphisms without any dominated splitting (recall example 1 in
section \ref{s5}).  
 The Main
Theorem can be re-casted in terms of the geometrical properties,
see Main Theorem Revisited in section \ref{ss17}.  In section \ref{ss121}, we adapt the Main
Theorem for the case that the endomorphism 
has partially
hyperbolic splitting, 
this is given in Theorem \ref{teo2} and the proof is an adaptation of \cite{PS1}.

In section \ref{s5}, we provide examples satisfying the main results.
Those satisfying the Main Theorem are done in such a way that they do not have any
dominated splitting and they are Derived from Expanding endomorphisms. For this case
we provide two type of examples: ones are built through bifurcation of periodic points
and the other are ``far from" expanding endomorphisms (see examples 1 and 2).
In examples 3 and 4, 
we show that there
are open sets of endomorphisms that satisfy Theorem \ref{teo2}. Those examples
are partially hyperbolic and they are not isotopic to expanding endomorphisms.


%




\section{Proof of the Main Result}\label{ss11}

Before starting the proof, we state a series of remarks that could help to understand the
hypotheses  of the Main Theorem and  in subsection \ref{ss12} we provide a sketch of the proof,
pointing out the main details and the general strategy.

%
%


\begin{obs}\label{obs7}
As we say in the introduction, the condition $diam(U_0)<1$ implies that we can assume that the closure of
$\widetilde{U}_0$ is contained in the interior of $[0,1]^n,$ where $\widetilde{U}_0$
is the lift of $U_0$ restricted to $[0,1]^n.$ Note that
$U_0$ do not need to be simply connected and could have
finitely many connected components. Actually the important fact is that
the closure of the convex hull of the lift
of $U_0$ restricted to $[0,1]^n$ is still contained in $(0,1)^n.$ Observe,
$diam_{int}(U_0^c)=diam_{int}(\mathfrak{U}_0^c),$ where
$\mathfrak{U}_0$ is the convex hull of $\widetilde{U}_0.$
\end{obs}


%



The Main Theorem is formulated for the $n-$dimensional torus. Some
of the examples provided in section \ref{s5} are isotopic to
expanding endomorphisms. Taking into account \cite{Shub}, we may
formulate the following:

\begin{conj}
The Main Theorem holds, at least, for any manifold supporting
expanding endomorphisms.
\end{conj}


\begin{obs}\label{obs2}
Using hypothesis (3) of the Main Theorem, 
given any point
$x\in U_1^c,$ we can construct a sequence $\{x_k\}_{k\geq 0}$
such that $x_0=x,$ $x_k\in U_1^c$ and $f(x_{k+1})=x_k$
for every $k\geq 0$. We call this sequence by \emph{inverse path}.
\end{obs}

\begin{obs}
The hypothesis of  diameter less than 1 and hypothesis (3) are technical. This means
that they are necessary conditions for the present proof of  our result,
 but we do not know if there exist weaker conditions that imply
 the thesis of our theorem.
\end{obs}



\begin{obs}\label{obs1} Observe that
$\Lambda_0:=\bigcap_{n\geq 0} f^{-n}(U_0^c)$
is an expanding  set. Moreover, from hypothesis (2) follows that  given any arc $\gamma$
in $U_0^c$ with diameter greater than $\delta_0$, there exists a point
$x\in\gamma$ such that  $f^k(x)$ is not in $U_1$ for any $k\geq 1.$
 Therefore, $\gamma\cap\Lambda_0\neq\emptyset$ and in particular $\Lambda_0$ is not trivial.
\end{obs}


\begin{defi}
Let $\Lambda$ be an expanding set  for $f\in E^1(M).$ 
If there is an open neighborhood $V$ of $\Lambda$ such that
$\Lambda=\bigcap_{k\geq 0} f^{-k}(\overline{V})$ then $\Lambda$
is said to be \emph{locally maximal} (or isolated) set.
$V$ is called the \emph{isolating block} of $\Lambda.$
\end{defi}

All previous remark can be summarized and extended in the next observation.

\begin{obs}\label{obs3}
Let us denote $\Lambda_1:=\bigcap_{n\geq 0} f^{-n}(U_1^c).$ This set has the
following pro\-per\-ties:
\begin{enumerate}
\item $\Lambda_1$ is an expanding  set.
\item By hypothesis (2) of the Main Theorem, 
        given any arc $\gamma$
       in $U_0^c$ with diameter greater than $\delta_0$, there exists a point
       $x\in\gamma$ such that  $f(x)\in\Lambda_1.$
\item Since the hypothesis $0<\delta_0<diam_{int}(U_0^c)$ is an open condition, we
       may take $U_1$ an open neighborhood of $\overline{U_0}$ such that
       $\delta_0<diam_{int}(U_1^c)<diam_{int}(U_0^c).$ Then for every arc $\gamma$ in $U_1^c$
       with diameter greater than $\delta_0$ holds that $\gamma\cap \Lambda_1$ is non empty.
\item $\Lambda_1$ is invariant, i.e. $f(\Lambda_1)=\Lambda_1$. It is clear that $\Lambda_1$ is
      forward invariant. So let us prove that $\Lambda_1\subset f(\Lambda_1).$
      Pick a point $x\in\Lambda_1$ and consider the sequence $\{x_k\}_{k\geq 0}$
      given by remark (\ref{obs2}). Let us show that $x_k\not\in W$ for any
      $k\geq 0,$ where
      $W=\cup_{n\geq 0} f^{-n}(U_1)=\Lambda_1^c.$ If this is not true,  there exist $k\geq 0$ and
      $n_k\geq 0$ such that $f^{n_k}(x_k)\in U_1.$ First, observe that
       remark (\ref{obs2}) implies that $f^n(x_k)=x_{k-n}$
      for $0\leq n\leq k.$ In particular, $f^k(x_k)=x_0$ if $k\geq 0.$
      And $f^n(x_k)=f^{n-k}(f^k(x_k))=f^{n-k}(x_0)$ for $n>k\geq 0.$
      Therefore,
      if $-k\leq -n_k\leq 0,$ then $f^{n_k}(x_k)=x_{k-n_k}.$
      Since every $x_k$ belongs to $U_1^c,$ we obtain that $f^{n_k}(x_k)$  belongs to $U_1^c$
      which is a contradiction because it was supposed that $f^{n_k}(x_k)\in U_1.$
      If  $-n_k< -k<0,$ then
        $f^{n_k}(x_k)=f^{n_k-k}(x_0).$ Since
      $x_0\in \Lambda_1,$ every positive iterate of $x_0$ by $f$ belongs to $U_1^c,$
      thus $f^{n_k}(x_k)\in U_1^c,$ which  contradicts the fact that $f^{n_k}(x_k)\in U_1.$
      Thus,  $x_k\in \Lambda_1$ for every $k\geq 0.$
\item In section \ref{ss13}, we prove that this set is locally maximal or it is contained in
      an expanding locally maximal set.
\end{enumerate}
\end{obs}


\subsection{Sketch of the Proof of  Main Theorem }\label{ss12}
Observe that if  $f$  satisfies the hypotheses 
of the Main Theorem,
then it satisfies the following property
denoted as \emph{internal radius
growing} (I.R.G.) property:

 \begin{quote}
\emph{There exists $R_0$
depending on the initial system such that given any
open set $U,$ there exist $x\in U$ 
and $K\in \mathbb{N}$ verifying  that $f^K(U)$ contains a ball of
a fixed radius $R_0$ centered in $f^K(x).$}
\end{quote}

In fact, since $f$ is volume expanding, then the lift of $f$
is a diffeomorphism in the universal covering space $\mathbb{R}^n$. In
consequence, given any open set $U\subset \mathbb{T}^n,$ volume
expanding implies that the diameter of the iterates by $f$ 
of $U$ grows on the covering space, 
(see Lemma \ref{afir4} for details). Then, for some $N>0,$
the diameter of $f^N(U)\cap U_0^c$ is greater or equal to
$\delta_0$ (the constant in the second hypothesis). Then we can pick an
arc in $f^N(U)\cap U_0^c$ of sufficiently large diameter
and using the second hypothesis we get
that there exists a point in $f^N(U)$ such that its forward orbit
remain in the expanding region. Therefore, the internal radius of
$f^{k+N}(U)$ grows as $k$ grows and the I.R.G. property follows.

Hence, if we have that $g$ also verifies the I.R.G. property, then the
Main Theorem is proved:  since every pre-orbit by $f$ is
dense in the manifold, given $0<\varepsilon<R_0,$ for $g$
$\varepsilon/2-$close to $f,$
 the pre-orbit of every point by $g$  are $\varepsilon-$dense (see subsection \ref{PMT}), then
$g^K(U)$ intersects $\{w\in g^{-n}(z):\, n\in\mathbb{N}\}$ for any $z$.
Therefore, taking pre-images by $g$, we get that $U$ intersects
$\{w\in g^{-n}(z):\, n\in\mathbb{N}\}$ for any $z$.

Therefore, the aim is to show that for every $g$ sufficiently close to $f$,
$g$ verifies the I.R.G. property, 
in other words we want to show that the I.R.G. property is robust.
In order
to prove this statement, we use a geometrical approach:
\begin{enumerate}
\item Since the initial map $f$ is volume expanding, then the perturbed map $g$
      is also volume expanding. So, its lift is also a diffeomorphism in the
      universal covering space $\mathbb{R}^n$.
\item Hypothesis (2) implies that there is an expanding
      subset $\Lambda_f$ that ``se\-pa\-ra\-tes", meaning that
      a nice class of arcs in $U_0^c$ intersect this set
      (see remarks \ref{obs1} and \ref{obs3} and lemma \ref{afir3} in section \ref{ss15}).
\item The set $\Lambda_f$ can be chosen as locally maximal
      (see lemma \ref{lema1} in section \ref{ss13}).
\item Hence the set $\Lambda_f$ has a continuation: for $g$ nearby, $f\mid_{\Lambda_f}$ is conjugate (see definition \ref{def2.5})
        to $g\mid_{\Lambda_g},$ and this conjugation is extended to a neighborhood of $\Lambda_f$
        and $\Lambda_g$
      (see propositions \ref{afir2} and \ref{extension} in section \ref{ss14}).
\item Therefore, the topological property of separation persists:
        for a nice class of arcs, every arc intersects $\Lambda_g$ following that the I.R.G. property holds for $g$
      (see lemma \ref{lema2} in section \ref{ss15}).
\end{enumerate}
\B

\subsection{Existence of an Expanding Locally Maximal Set for $f$}\label{ss13}

In the present subsection (lemma \ref{lema1}) we show that $\Lambda_1$ (as defined in remark \ref{obs3}) is either locally maximal or is contained in a locally maximal one.
The third hypothesis in the Main Theorem is essential to prove this fact (see remark \ref{cro-fish} for a discussion about this issue).

\begin{lema}\label{lema1}
Either $\Lambda_1$ is a locally maximal set or there exists
$\Lambda^*$ an expanding  locally maximal set for $f$ such that
$\Lambda_1\subset\Lambda^*$ and
$\Lambda^*$ verifies that  every arc $\gamma$
in $U_0^c$ with diameter greater than $\delta_0$ has a point
such that the image by $f$ belongs to $\Lambda^*.$ Moreover,
every arc $\gamma$
in $U_1^c$ with diameter greater than $\delta_0$ intersects $\Lambda^*.$
\end{lema}

\pp
We may divide the proof in two cases:

\textbf{Case I.} $\Lambda_1\cap \partial U_1=\emptyset.$

Let us observe that $\Lambda_1\cap \partial U_1=\emptyset$ implies
that $\Lambda_1$ is contained in the open neighborhood $V=int(U_1^c).$
Then $V$ is an isolating block
for $\Lambda_1,$ therefore $\Lambda_1$ is locally maximal.

\textbf{Case II.} $\Lambda_1\cap \partial U_1\neq\emptyset.$

In this case, $V$ fails to be an isolating neighborhood.
To overcome this situation, we extend $V$ in a proper way and we show that the
extension is now an isolating neighborhood of the respective maximal invariant set.
Choose $\varepsilon>0$ sufficiently small such that the open ball
 $\mathbb{B}_{\varepsilon}(x)$ is
contained in  $U_0^c$ for
all $x\in \Lambda_1$ and  for every $x\in \Lambda_1,$ since $f$ is a local diffeomorphism,
there exists an open set $U_x$ such that
$f\mid_{U_x}:U_x\rightarrow \mathbb{B}_{\varepsilon}(x)$ is a diffeomorphism.
 Note that
the collection $\{\mathbb{B}_{\varepsilon}(x)\}_{x\in \Lambda_1}$ is an open cover
of $\Lambda_1.$
Since $\Lambda_1$ is compact, there is a finite subcover, let us say
$\{\mathbb{B}_{\varepsilon}(x_i)\}_{i=1}^{N}.$

Fix $\lambda_0^{-1}<\lambda'<1,$ where $\lambda_0$ is the
expansion constant of $f$ and pick $N'$ greater or equal to $N,$
the cardinal of the finite subcover of $\Lambda_1,$ such that
for every $y\in\Lambda_1$, there is  $i=i(y)\in\{1,\ldots, N'\}$ such that
$\mathbb{B}_{\lambda'\varepsilon}(y) \Subset \mathbb{B}_{\varepsilon}(x_i),$ i.e.
$\overline{\mathbb{B}_{\lambda'\varepsilon}(y)}\subset \mathbb{B}_{\varepsilon}(x_i).$

Let us define $W=\displaystyle\bigcup_{i=1}^{N'}\mathbb{B}_{\varepsilon}(x_i)$ and
$\widehat{W}=\displaystyle\bigcup_{i=1}^{N'}\overline{\mathbb{B}_{\varepsilon}(x_i)}.$

By remark (\ref{obs3}) $\Lambda_1$ is invariant, then
 we have that for every $x_i,$ there exists at least one
$x_i^j\in \Lambda_1$ such that $f(x_i^j)=x_i.$ Let us consider for every $1\leq i\leq N'$
all the possible pre-images by
$f$ of $x_i$ that belongs to $\Lambda_1,$ i.e. recall that  $f$ is a local diffeomorphism,
hence for every point $x\in M,$ the cardinal $\sharp\{f^{-1}(x)\}=N_f$ is constant, then
for every $i\in \{1,\ldots,N'\},$ there exist $K_i\subset \{1,\ldots,N_f\}$
such that if $j\in K_i$ then $x_i^j\in \Lambda_1$  and $f(x_i^j)=x_i.$
Therefore for every $i\in\{1,\ldots, N'\}$ and
for every $j\in K_i,$
there exist open sets $U_i^j$ such that $x_i^j\in U_i^j$ and
$f\mid_{U_i^j}: U_i^j \rightarrow \mathbb{B}_{\varepsilon}(x_i)$ is a diffeomorphism.
Given $i\in\{1,\ldots, N'\},$  for every $j\in K_i$ consider the inverse branches,
$\varphi_{i,j}:\mathbb{B}_{\varepsilon}(x_i)\rightarrow U_i^j $
such that
$$
\begin{array}{rl}
  \varphi_{i,j}(x_i)= & x_i^j,\\
  f \circ \varphi_{i,j}(x)=& x,\quad \forall \,x\in \mathbb{B}_{\varepsilon}(x_i).
\end{array}
$$

 Now, consider
 $\Lambda^*=\bigcap_{n\geq 0} f^{-n}(\widehat{W}).$
 Clearly, $\Lambda_1\subset \Lambda^*\subset U_0^c$ and $\Lambda^*$ is an expanding set.
 In order to show that $\Lambda^*$ is locally maximal, it is enough to
 show that $\Lambda^*\cap \partial \widehat{W}=\emptyset,$
 which is equivalent showing that $f^{-1}(\widehat{W})$ is  contained in $W.$
 Just to make more clear what follows, let us rewrite $f^{-1}(\widehat{W})$  in terms of the
 inverse branches,
 $$f^{-1}(\widehat{W})=f^{-1}(\bigcup_{i=1}^{N'}\overline{\mathbb{B}_{\varepsilon}(x_i)})
 = \bigcup_{i=1}^{N'}\bigcup_{j\in K_i}\overline{ \varphi_{i,j}(\mathbb{B}_{\varepsilon}(x_i))}.$$
 So, it is enough to show that
 $$\overline{ \varphi_{i,j}(\mathbb{B}_{\varepsilon}(x_i))}\subset \mathbb{B}_{\varepsilon}(x_{m_{i,j}}),$$
 for some $x_{m_{i,j}}\in \{x_1,\ldots,x_{N'}\}.$
In fact,
$$\varphi_{i,j}(\mathbb{B}_{\varepsilon}(x_i))=U_i^j\subset \mathbb{B}_{\lambda_0^{-1}\varepsilon}(\varphi_{i,j}(x_i))
\subset \mathbb{B}_{\lambda'\varepsilon}(\varphi_{i,j}(x_i))=\mathbb{B}_{\lambda'\varepsilon}(x_i^j)$$
then, there exists $m_{i,j}\in\{1,\ldots,N'\}$ such that
$\mathbb{B}_{\lambda'\varepsilon}(x_i^j)\Subset \mathbb{B}_{\varepsilon}(x_{m_{i,j}}),$
and the assertion holds.

Since  $\Lambda_1\subset \Lambda^*\subset U_0^c$ and $\Lambda_1$
intersects the image by $f$ of every arc
$\gamma$ in $U_0^c$ with diameter larger than $\delta_0,$ then
follows that $\Lambda^*$ also verifies the latter property. In
particular, $\Lambda^*$ intersects every arc $\gamma$ in $U_1^c$
with diameter larger than $\delta_0.$
\B


\begin{figure}[h]
\begin{center}
\includegraphics[scale=0.3]{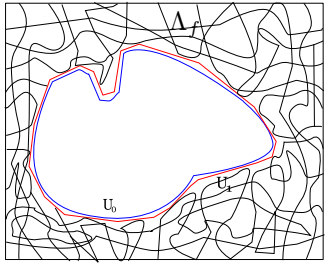}
\rput[0,0](-2.3,3){ $\Lambda_f$}
\end{center}
\caption{\textit{$\Lambda_f$ looks like a net which is an expanding set that ``separates''}}\label{graf1.3}
\end{figure}

\begin{obs}\label{cro-fish}
We want to highlight that for diffeomorphisms
there exist examples of hyperbolic sets that are not contained
in any locally maximal hyperbolic set, see for instance \cite{Crovisier}
and \cite{Fisher}. A similar construction seems feasible for endomorphisms.
In our context,  hypothesis (3) allows to overpass this problem and   guarantees
 that $\Lambda_1$ is an invariant set. Moreover, we can consider a finite covering
 $\{\mathbb{B}_{\varepsilon}(x_i)\}_{i=1}^{N'}$ for
 $\Lambda_1,$ with $x_i\in \Lambda_1,$ in such a way that for every point $y\in\Lambda_1,$ there is  $x_i$
 such that  $\overline{\mathbb{B}_{\lambda'\varepsilon}(y)}\subset \mathbb{B}_{\varepsilon}(x_i).$
 Thus we conclude that $\Lambda^*$ is contained
 in the interior of $\widehat{W}$ and therefore
  the expanding set $\Lambda_1$ is either locally maximal or is contained in a locally maximal
 expanding set.
\end{obs}


\subsection{Continuation of the Expanding Locally Maximal Set}\label{ss14}

First, in proposition \ref{afir2}  
we prove that
$g\mid_{\Lambda_g}$ is conjugate to $f\mid_{\Lambda_f},$ where
$\Lambda_g$ is the maximal invariant set  associated to $g,$ for $g$
sufficiently close to $f.$
 This is standard in hyperbolic theory using a
 shadowing's lemma argument. We provide  the proof for completeness and to show how
 the conjugacy can be extended to a neighborhood. This is done in proposition \ref{extension}.
 We want to remark that to construct the topological conjugacy between
 $g\mid_{\Lambda_g}$ and $f\mid_{\Lambda_f}$ is not necessary
 that $\Lambda_f$ be locally maximal invariant, however, this property is essential in the proof of proposition \ref{extension}.
\begin{defi}
The sequence $\{x_n\}_{_{n\in\mathbb{Z}}}$ is said to be a $\delta-$\emph{pseudo orbit}
for $f$ if $d(f(x_n), x_{n+1})\leq \delta$ for every $n\in\mathbb{Z}.$
\end{defi}

\begin{defi}
We say that a $\delta-$pseudo orbit $\{x_n\}_{_{n\in\mathbb{Z}}}$ for $f$
is $\varepsilon-$\emph{shadowed} by a full orbit $\{y_n\}_{_{n\in\mathbb{Z}}}$ for $f$
if $d(y_n, x_n)\leq \varepsilon$ for every $n\in\mathbb{Z}.$
\end{defi}

It follows that the \emph{Shadowing Lemma} holds for $C^1$ endomorphisms. 

\begin{lema}\label{sl}
Let $M$ be a Riemannian manifold, $U\subset M$ open, $f:U\rightarrow M$ a $C^1$ expanding endomorphism, and
$\Lambda\subset U$ a compact  invariant expanding set for $f$. Then there exists a neighborhood
$\mathcal{U}(\Lambda)\supset \Lambda$
such that whenever $\eta>0$ there is an $\varepsilon>0$ so that every $\varepsilon-$pseudo orbit for $f$
in $\mathcal{U}(\Lambda)$ is $\eta-$shadowed by a full orbit of $f$ in $\Lambda.$
If $\Lambda$ is locally maximal invariant set,
then the shadowing full orbit is contained in $\Lambda$.
\end{lema}

\pp
For details, see for instance \cite{Liu}.

\begin{defi}\label{def2.5}
Let $f:M\rightarrow M$ and $g:N\rightarrow N$ be two maps and let
$\Lambda_f$ and $\Lambda_g$ be invariant sets by $f$ and $g$
respectively. We say that $f:\Lambda_f \rightarrow \Lambda_f$ is
\emph{topologically conjugate} to $g:\Lambda_g \rightarrow
\Lambda_g$ if there exists a homeomorphism (in the relative
topology) $h:\Lambda_f \rightarrow \Lambda_g$ such that $h\circ f=
g\circ h.$
\end{defi}

This a typical notion for hyperbolic sets, see \cite{Shub1}.


In order to fix some notation in what follows,
we denote by
 $\Lambda_f$ the expanding  locally maximal set for $f:$
  $\Lambda_f$ is either $\Lambda_1,$ in the case it is locally maximal, or
 it is  $\Lambda^*$ given in Lemma \ref{lema1}. We also denote by $U$ the isolating
  block of   $\Lambda_f.$

\nocite{LC}

\begin{prop}\label{afir2}
There exists $\mathcal{V}_1(f)$
an open neighborhood of $f$ in $E^1(\mathbb{T}^n)$ such that if $g\in \mathcal{V}_1(f),$ then $g$
is expanding on
$\Lambda_g=\bigcap_{n\geq 0} g^{-n}(U)$
and there exists an homeomorphism $h_g:\Lambda_g \rightarrow \Lambda_f$ 
that conjugate $f\mid_ {\Lambda_f}$ and $g\mid_ {\Lambda_g}$ and $h_g$  is close to the identity.
\end{prop}

\pp
In order to get the conjugacy we use the
Shadowing Lemma for $C^1$ expanding endomorphisms, see lemma \ref{sl}.

Since $\Lambda_f$ is an expanding  locally maximal set for $f,$
there exists 
$\beta>0$
such that 
$f$ is expansive with constant $\beta$
in $\Lambda_f.$

Fix $0<\eta<\beta.$ By the endomorphism version of the Shadowing Lemma,
there exists $\varepsilon>0$
such that 
any $\varepsilon-$pseudo orbit for $f$ within
$\varepsilon$ of $\Lambda_f$
is uniquely $\eta-$shadowed by a full orbit in $\Lambda_f.$

Take $N$ such that
$$\bigcap_{j=0}^N f^{-j}(U)\subset \{q:\,d(q,\Lambda_f)<\varepsilon/2\}.$$
There exists a $C^0$ neighborhood $\mathcal{V}(f)$ of $f$ such that for $g$ in $\mathcal{V}(f)$
$$\bigcap_{j=0}^N g^{-j}(U)\subset \{q:\,d(q,\Lambda_f)<\varepsilon/2\}$$
and
for any  $x\in\bigcap_{j=0}^N g^{-j}(U),$
we may consider $\{x_n\}_{_{n\in\mathbb{Z}}}$ a full orbit for $g,$ where $x_0=x,$ getting that
$\{x_n\}_n$ is an $\varepsilon-$pseudo orbit for $f.$

Let $\Lambda_g=\bigcap_{n\geq 0} g^{-n}(U).$ Taking an open subset $\mathcal{V}_1(f)$ of $\mathcal{V}(f)$
small enough in the $C^1$ topology, then for $g\in\mathcal{V}_1(f),$ $\Lambda_g$ is an expanding locally maximal
set for $g$. If $g$ is close enough to $f$, then $g$ is also expansive with constant $\beta.$ Moreover,
the Shadowing Lemma also holds for $g.$

Take $g\in \mathcal{V}_1(f).$ Given $x\in\Lambda_g,$
consider $\{x_n\}_{_{n\in\mathbb{Z}}}$
a full orbit for $g,$ where $x_0=x.$ As  $\{x_n\}_n$ is an $\varepsilon-$pseudo orbit for $f,$
there exists a unique full orbit $\{y_n\}_{_{n\in\mathbb{Z}}}$ for $f$ with $y_0=y\in\Lambda_f$ that
$\eta-$shadows $\{x_n\}_{_{n\in\mathbb{Z}}}$.


Let us define $h_g:\Lambda_g\rightarrow\Lambda_f$ by $h_g(x)=y,$ where $y$ is given by the Shadowing Lemma.
 By the uniqueness of the
shadowing point, this map is well defined. The continuity of $h_g$ follows from the shadowing lemma.

Moreover, $h_g\circ g=f\circ h_g.$ In fact, consider the sequence $\{z_n\}_{_{n\in\mathbb{Z}}}$ where
$z_n=g(x_n)=x_{n+1}.$ This $\varepsilon-$pseudo-orbit is $\eta-$shadowed by a unique full orbit
$\{w_n\}_{_{n\in\mathbb{Z}}}$ for $f,$ with $w_0=w\in\Lambda_f.$ Then,  for every $n\in\mathbb{Z},$
$$
\begin{array}{ll}
 d(w_n,z_n)&= d(f^n(w_0),x_{n+1}) = d(f^n(h_g(z_0)),g^n(g(x_0))) \\
           &= d(f^n(h_g\circ g(x_0)),g^n(g(x_0)))\\
& = d(f^{n+1}\!\!\circ\! f^{-1}\!\circ \!h_g\!\circ\! g(x_0),g^{n+1}(x_0))<\eta
\end{array}
$$
Observe that $f^{-1}\circ h_g\circ g(x_0)=w_{-1}$ is $\eta-$shadowing $x_0.$ So, by uniqueness, we have
that $f^{-1}\circ h_g\circ g(x_0)=y_0;$ i.e. $h_g\circ g(x)=f\circ h_g(x).$

Since  we can apply the Shadowing
Lemma  for $\Lambda_g$ using the same constants as in the construction of $h_g$, we define a map
$h_f:\Lambda_f\rightarrow \Lambda_g$ such that $h_f\circ f=g\circ h_f.$
In fact, if $\{y_n\}_{n\in\mathbb{Z}}$ is a full orbit for $f$ with $y_0\in \Lambda_f,$
then it is an $\varepsilon-$pseudo orbit for $g.$ Hence, this pseudo orbit is
uniquely shadowed by a full orbit $\{x_n\}_{n\in\mathbb{Z}}$ for $g,$ with $x_0\in \Lambda_g.$
Thus, $h_f(y_0)=x_0$ and $d(y_n,x_n)<\eta$ for every $n\in\mathbb{Z};$
moreover, $h_f$ is continuous  and satisfies
$h_f\circ f=g\circ h_f$ just as $h_g.$

Next, let us verify that $h_g$ is one to one. Let $p_1,\;p_2\in \Lambda_g$ be two points such
that $h_g(p_1)=h_g(p_2).$ Note that
$d(f^n(h_g(p_1)),g^n(p_1))<\eta$ and $d(f^n(h_g(p_2)),g^n(p_2))<\eta$ by construction.
Then $h_g(p_1)$ is $\eta-$shadowed by $p_1$ and $p_2,$ 
which by uniqueness gives that $p_1=p_2.$

Finally, for $y\in \Lambda_f,$ consider a full orbit of $h_f(y)$ by $g.$
Since $d(g^n(h_f(y)),f^n(y))$ is small for all $n$ and some
$f$ full orbit of $y$ shadows the $g$ full orbit of $h_f(y),$ we have
that $h_g(h_f(y))=y.$ Hence, $h_g$ is onto and therefore is a homeomorphism.
\B

The next proposition  is a version for expanding endomorphisms of a result  already provided
for the case of hyperbolic diffeomorphisms in \cite[Theorem 4.1]{Robinson}.
The goal is to show that we can extend the conjugation between $f\!\!\mid_{\Lambda_f}$
and $g\!\!\mid_{\Lambda_g}$ to an open neighborhood $U$ of $\Lambda_f$ in such a way that
still is a homeomorphism that conjugate $f\!\!\mid_U$ and $g\!\!\mid_U,$ noting that
the conjugacy is unique just in $\Lambda_f.$  We are going
to use this extension in next section for proving that the property of $\Lambda_f$ disconnecting
a ``nice'' class of arcs is robust.

\begin{prop}\label{extension}
The homeomorphism $h_f:\Lambda_f\rightarrow \Lambda_g$ in proposition
(\ref{afir2}) can be extended as a homeomorphism $H$
to an open neighborhood of $\Lambda_f$ such that $H\circ f=g\circ H.$
\end{prop}

\pp
This  proof is inspired in the geometrical approach  given by Palis in \cite{Palis}
and also used to prove the Grobman-Hartman Theorem in \cite[pp.96]{Shub1}.

Other alternative proof consist in using inverse limit space, in such a way that the expanding set $\Lambda_f$
becomes a hyperbolic set for a diffeomorphism and  so Theorem 4.1 in \cite{Robinson} can be applied.
Observe, that to have a well defined inverse limit and that the induced set associated to $\Lambda_f$
verifies the hypothesis of the mentioned theorem in \cite{Robinson} is needed that $\Lambda_f$ is locally maximal invariant.

The goal is to choose  an appropriate  isolating neighborhood  $U$ of $\Lambda_f$ and to construct an homeomorphism
from $U$ onto itself, using the inverse branches of $f$ and $g$, first defined in   a fundamental domain $D_f$
 for $f$ (i.e. a set $D_f$ such that for every $x\in U\setminus\Lambda_f,$ there exists $n\in\mathbb{N}$ such that
$f^n(x)\in D_f$) and then extended to $U$ using inverse iterates. 
Observe that the isolating block of $\Lambda_f$ is also an isolating block of $\Lambda_g.$
Now we can take a fundamental domain for $g,$ $D_g,$ as it was done for $f.$
Note that $D_f$ is defined as $U\setminus f^{-1}(U)$ 
and since $f^{-1}(U)$ is properly contained in $U,$ it follows that the same holds for
$g$ and therefore $D_f$ and $D_g$ are homeomorphic.
Then it is taken an homeomorphism $H$
between both fundamental domains $D_f$ and $D_g.$ This homeomorphism is saturated to $U\setminus\Lambda_f$
by backward iteration, i.e. if $x\in U\setminus \Lambda_f,$
let $n$ be such that $f^n(x)\in D_f,$ take $H\circ f^n(x)$ and then $g^{-n}\circ H\circ f^n(x)$ where
$g^{-n}$ is taken carefully using the corresponding inverse branches.

 Denote by
$N_f$ the cardinal of $\{w\in f^{-1}(x)\}$. Since $f$ is a local
diffeomorphism, $N_f$ is cons\-tant. Let $K\subset \{1,\ldots,
N_f\}$ be such that for every $i\in K,$ there exist $U_i^f\subset
U$ and $\varphi_i^f:U\rightarrow U_i^f$ inverse branch of $f$ such
that $\varphi_i^f(U)=U_i^f$ and $f(U_i^f)=f\circ\varphi_i^f(U)=U.$
Also, for $g$ as in proposition \ref{afir2}, for every $i\in K,$
there exist $U_i^g\subset U$ and $\varphi_i^g:U\rightarrow U_i^g$
the inverse branch of $g$  such that $\varphi_i^g(U)=U_i^g$ and
$g(U_i^g)=g\circ\varphi_i^g(U)=U.$

To construct an homeomorphism $H$ on $U$  satisfying
$H\circ f= g\circ H$ and $H\mid_{\Lambda_f}=h_f$
we can begin as follows.  Suppose that
the restriction
$H: \partial U\rightarrow \partial U$ is any well-defined
orientation preserving diffeomorphism.
The restriction of $H$ to $\partial U_i^f$ is then defined
as follows $H(x)= \varphi_i^g\circ H \circ f(x)$
if $x\in \partial U_i^f$ because $H$ conjugate $f$ and $g$.
Now we extend $H$ to a diffeomorphism which send $U\setminus \bigcup_{i\in K}U_i^f$
bounded by $\partial U$ and $\partial U_i^f$ onto
$U\setminus \bigcup_{i\in K}U_i^g$
bounded by $\partial U$ and $\partial U_i^g.$
Since we may assume that the Hausdorff distance between $U$
and $\Lambda_f$ is small, see lemma \ref{lema1}, then the initial $H$
is close to the identity. Let us say that $d(H(x),x)<\eta,$
where $\eta>0$ is given arbitrarily.

Given $i,j\in K,$ denote $U_{j,i}^f=\varphi_j^f\circ \varphi_i^f(U)$ and
$U_{2\;i}^f=U_i^f\setminus \bigcup_{j\in K} U_{j,i}^f.$
If $x\in \partial U_{j,i}^f$ then
$H(x)=\varphi_j^g\circ \varphi_i^g\circ H\circ f^2(x)\in \partial U_{j,i}^g.$
If $x\in U_{2\;i}^f\setminus \Lambda_f$ then $H(x)=\varphi_i^g \circ H\circ f(x)\in U_{2\;i}^g.$

Doing this process inductively we have that:
Given $i_1,\ldots,i_n\in K,$ denote $U_{i_n,\ldots,i_1}^f=\varphi_{i_n}^f\circ\cdots\circ \varphi_{i_1}^f(U)$ and
$U_{n\;(i_{n-1},\ldots,i_1)}^f=U_{i_{n-1},\ldots,i_1}^f\setminus \bigcup_{i_n\in K} U_{i_n,\ldots,i_1}^f.$
If $x\in \partial U_{i_n,\ldots,i_1}^f$ then
$H(x)=\varphi_{i_n}^g\circ \cdots\circ \varphi_{i_1}^g\circ H\circ f^n(x).$
If $x\in U_{n\;(i_{n-1},\ldots,i_1)}^f\setminus \Lambda_f$ then
$H(x)=\varphi_{i_{n-1}}^g\circ \cdots\circ \varphi_{i_1}^g \circ H\circ f^{n-1}(x).$
And $H(x)=h_f(x)$ if $x\in\Lambda_f.$

Let us prove that $H$ is continuous.

Given $x\in\Lambda_f,$ let $(x_n)_n$ be a sequence in $U\setminus \Lambda_f$ such
that $x_n\rightarrow x,$ when $n\rightarrow \infty.$
Let us prove that $H(x_n)\rightarrow H(x),$ when $n\rightarrow \infty.$

First, consider $\{z_k\}_{k\in\mathbb{Z}}$ an $f-$full orbit in $\Lambda_f$ such that
$z_0=x$ and for every $n\in\mathbb{N},$ consider
 $\{z_k^n\}_{k\in\mathbb{Z}}$ a full orbit by $f$ associated to each $x_n$ using the corresponding
 inverse branches (for the backward iterates) given by the full orbit of $x,$ where
$z_0^n=x_n.$
Since $f$ is continuous, for every $k\in\mathbb{Z},$  we have that $z_k^n\rightarrow z_k$
when $n\rightarrow \infty.$

Note that for every $n\in\mathbb{N},$ there exists $k_n>0$ such that
$z_{k_n}^n\in U\setminus \bigcup_{i\in K}U_i^f.$ Furthermore,
$z_k^n\in U$
for every $k\in [-k_n,k_n].$
Since $H\circ f=g\circ H,$ we get that $H(x_n)\in \bigcap_{k=-k_n}^{k_n}g^k(U).$

Hence, for $\eta$ and $\varepsilon$  as in proposition \ref{afir2} and for every $n\in\mathbb{N},$ we have
that $\{z_k^n\}_{k=-k_n}^{k_n}$ is a finite $\varepsilon-$pseudo orbit for $g$
and it is $\eta-$shadowed by a $g-$orbit of $H(x_n)$ until $k_n$ for forward iterates and
$-k_n$ for backward iterates.

Observe that as $m$ goes to infinity, the finite pseudo orbit $y_n^m=\{z_k^n\}_{k=-m}^m$
becomes longer.
Consider now the sequence $\{y_n^m\}_n.$
Then $y_n^m\rightarrow \{z_k\}_{k=-m}^m$ when $n\rightarrow\infty.$ Hence,
the sets of shadowing points of the finite pseudo orbits $y_n^{k_n}$
converge to the shadowing point of the
infinite pseudo orbit $\{z_k\}_k$, then $H(x_n)\rightarrow h_f(x)=H(x)$ when
$n\rightarrow \infty.$ \B


\subsection{The Locally Maximal Set ``Separates"}\label{ss15}

The main goal of this section is to show that the locally maximal set for $f$
has a topological property that persists under perturbation, roughly speaking
means that $\Lambda_f$ and $\Lambda_g$ disconnect small open sets. Using that we prove that $\Lambda_f$
intersects  ``some nice" class of arcs in $U_1^c$ and which  also intersect $\Lambda_g$
 for all $g$ nearby $f$. The first question that arise is: which arcs belong to this
``nice" class? The second questions in the context of proving the Main Theorem is: why is this property enough?
The third question is: why does the ``nice class" of arc exist?
 All these questions are answered along the section, but to give some
brief insight about the main ideas we want to make some comments:

\begin{enumerate}
\item Roughly speaking, these ``nice arcs'' are arcs that have the property that can be used
     to build ``nice cylinders" (see definition \ref{def1.18})
    containing the initial arc and such that $\Lambda_f$ ``separates" (see definition \ref{def1.20})
        this cylinder in a ``robust way" (see lemmas \ref{afir3} and \ref{lema2}).
\item It is enough to consider these ``nice" class of arcs to finish the proof of the Main Theorem.
    Suppose that the existence  of this class of arcs is proved and suppose that given any open set,
    there is an iterate by $g$ that contains a ``nice" arc (see claim \ref{afir6} and lemma \ref{afir4}). Then  there is a point
in this iterate which its forward orbits stay in the expanding
region and arguing as in the beginning of subsection \ref{ss12}
 it is concluded the density of the pre-orbit of any point for the perturbed map.
\item Therefore, to finish, we show in claim (\ref{afir6}) that every large arc admits a ``nice" arc. Later it is shown that any open set has an iterate, in the universal covering,  with an arbitrary large arc (see lemma  \ref{afir4}).
\end{enumerate}

Let us define the concepts involved in this section.

\begin{defi}\label{cyl}(\textbf{Cylinder})
Given a differentiable arc  $\gamma$ and $r>0,$ it is said that
$C(\gamma,r)$ is a \emph{cylinder} centered at $\gamma$ with
radius $r$ if $$C(\gamma,r):=\bigcup_{x\in\gamma}
([T_x\gamma]^{\perp})_r,$$ where  $([T_x\gamma]^{\perp})_r$
denotes the closed ball 
centered at $x$ with
radius $r$ intersected with
$[T_x\gamma]^{\perp}$ the orthogonal to the tangent to $\gamma$ in
$x.$
\end{defi}

\begin{defi}(\textbf{Simply connected cylinder})
Given a differentiable arc  $\gamma$ and $r>0,$
it is said that a cylinder $C(\gamma,r)$ is \emph{simply connected}
if  it is  retractile to a point.
\end{defi}

\begin{obs}
Fixed the radius, a cylinder as defined in \ref{cyl} could be not retractile to a point.
In this case, working in the universal covering
space, consider the convex hull of its lift and then project it on the
manifold. We
call the resulting set as \emph{simply connected cylinder} as well and is denoted
in the same way as above.
\end{obs}

\begin{figure}[h]
\begin{center}
\includegraphics[scale=0.2]{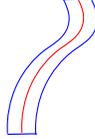}
\end{center}
\vspace*{-4mm}
\caption{\textit{Nice cylinder}}\label{graf1.1a}
\end{figure}

\begin{defi}(\textbf{Nice cylinder})\label{def1.18}
Given  an arc $\gamma$ and $r>0$, it is said that a cylinder
$C(\gamma,r)$ is a \emph{nice cylinder} if it is simply connected cylinder
and  if $x_A$ and $x_B$ are the extremal points of $\gamma$ then
$A:=([T_{x_A}\gamma]^{\perp})_r\subset \partial C(\gamma,r)$
and $B:=([T_{x_B}\gamma]^{\perp})_r\subset \partial C(\gamma,r).$
In this case, we say that $A$ and $B$ are the \emph{top and bottom sides}
of the cylinder. (See figure \ref{graf1.1a}).
\end{defi}

\vspace*{-8mm}

\vspace*{-3mm}
\begin{obs}
In general, given any cylinder, as defined in \ref{cyl}, it has
not necessarily  top and bottom sides, and may 
not be 
simply connected.
\end{obs}

Hereafter,
fix $U_2$ an open set such that $\overline{U_1} \subset U_2,$ where 
$U_1$ is the same as in hypothesis (2) in the Main Theorem,
and $\delta_0<diam_{int}(U_2^c)<diam_{int}(U_0^c).$ Let
$d_1=d_H(U_2, U_1)>0,$ where $d_H$ denotes the
Hausdorff metric, and let $k\in \mathbb{N}$ such that
$\delta'_0=\delta_0+\frac{d_1}{3k}<diam_{int}(U_2^c).$

Let us denote by $\widetilde{U}_0$ the lift of  $U_0,$
$\pi$ the projection of $\mathbb{R}^n$ onto $\mathbb{T}^n$ and
$\mathfrak{U}_0$ the convex hull of $\widetilde{U}_0\cap [0,1]^n$.
 Consider $\mathrm{P}_i(\mathfrak{U}_0)$
the projection of $\mathfrak{U}_0$ in the $i-$th coordinate in the
$n-$dimensional cube $[0,1]^n.$ Since $diam(U_0)<1$ and remark (\ref{obs7}),
for every $1\leq i\leq n,$ there exist $0<k_i^-<k_i^+<1$ such that
$k_i^-<\mathrm{P}_i(\mathfrak{U}_0)<k_i^+.$
Note that $1+k_i^--k_i^+>\delta'_0$
for every $i$, because $1+k_i^--k_i^+> diam_{int}(U_0^c)$ by construction.

Let $R_i^m=\{x\in\mathbb{R}^n: k_i^-+m<x_i<k_i^++m\}$ with $m\in\mathbb{Z},$ $1\leq i\leq n$
and $x_i$ is the $i-$th coordinate of $x$. Thus,
$\mathfrak{U}_0\subset\displaystyle\bigcap_{m\in\mathbb{Z}, 1\leq i\leq n} R_i^m.$
Denote by $L_i^+=\{x\in\mathbb{R}^n: x_i=k_i^+\}$ and $L_i^-=\{x\in\mathbb{R}^n: x_i=k_i^-\}.$
Let $\tilde{f}$ be the lift of $f$.

%


The next claim answer the third question stated at the beginning of the section.

\begin{figure}[h]
\begin{center}
\includegraphics[scale=0.34]{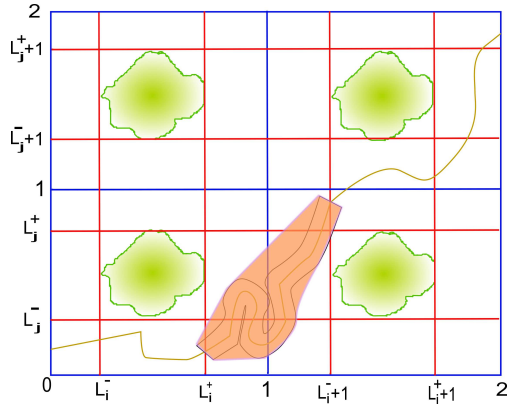}
\end{center}
\vspace*{-5mm}
\caption{\textit{Every arc admits a sub-arc with a nice cylinder}}\label{graf1.1}
\end{figure}

\begin{afir}\label{afir6}
Let  $m>2\sqrt{n}$ be fixed.
Given any arc $\gamma$ in $\mathbb{R}^n$ with $diam(\gamma)>m,$
 there exist an arc
$\gamma'\subset \gamma,$  $1\leq i\leq n$ and $j\in\mathbb{Z}$
such that $\partial \gamma'\cap (L_i^++j)\neq \emptyset,$ $\partial \gamma'\cap
(L_i^-+j+1)\neq \emptyset$ and $P_i^j(\gamma')\subset [k_i^++j, k_i^-+j+1],$
where $P_i^j(\gamma')$ denotes the projection of $\gamma'$ on the
interval $[j,j+2]$ of the $i-$th coordinate. 
Moreover, $\gamma'$ admits a nice cylinder $\gamma^*=\pi(\gamma')$
in  $U_2^c,$ with diameter of $\gamma^*$ larger than $\delta_0$
and also admitting  a nice cylinder contained in $U_1^c$. (See
figure \ref{graf1.1})
\end{afir}

%

\pp
Take $\gamma$ an arc  with diameter larger than $m,$
then
the projection of $\gamma$ in the $i-$th coordinate  contains an  interval
of the kind formed by $k_i^+$ and
$1+k_i^-$ for some $i$ (or formed by $k_i^++j$ and $k_i^-+j+1$ for some $j\in\mathbb{Z}$).
If it is not true, $\gamma$
would be in a $n-$dimensional cube with sides smaller than $k_i^+-k_i^-<1$
and this cube has diameter smaller than $\sqrt{n},$ but this contradict the
fact that $diam(\gamma)>m.$ Hence, we may pick an arc $\gamma'$ in $\gamma$
such that
$\partial \gamma'\cap (L_i^++j) \neq \emptyset,$  $\partial \gamma'\cap (L_i^-+j+1) \neq \emptyset$
and $P_i^j(\gamma')\subset [k_i^++j, k_i^-+j+1]$ for some $1\leq i\leq n$ and some $j\in\mathbb{Z}$.
Therefore, diameter of $\gamma'$ is greater than $\delta_0$ and in consequence
its projection in $\mathbb{T}^n$ also
has diameter greater than $\delta_0.$

Moreover, since the projection of
$\gamma'$ by $\mathrm{P}_i$ %
is in between
$k_i^++j$ and $k_i^-+1+j,$ we may construct a   cylinder centered at $\gamma'$
and radius $\frac{d_1}{2}$ such that this cylinder is ``far" away from $\widetilde{U}_0,$
so this cylinder could be simply connected or, if it is not simply connected cylinder,
it has  holes that are different from
$\widetilde{U}_0.$ In the case that the cylinder is not simply connected, we
consider the convex hull of the cylinder, since the original cylinder is bounded by $L_i^++j$ and
$L_i^-+j+1,$ then the convex hull stay in between these two hyperplanes and therefore it does not intersect
$\widetilde{U}_0.$
By abuse of notation, let us denote this set by $C(\gamma',\frac{d_1}{2}),$ it is a simply connected cylinder.
 Observe that by construction, this
cylinder will have top and bottom sides, thus $C(\gamma',\frac{d_1}{2})$ is a nice cylinder.

Take $\gamma^*=\pi(\gamma'),$ note that $\gamma'$ can be choose such that $\gamma^*$ is contained in  $U_2^c$
and the diameter of $\gamma^*$ is larger than $\delta_0,$
 then projecting the nice cylinder for $\gamma'$ in $\mathbb{T}^n$  we obtain a nice cylinder for  $\gamma^*$
 which is denoted by $C(\gamma^*,\frac{d_1}{2}).$ This nice cylinder has the property that every arc that
 goes from  bottom to top side has diameter at least $\delta_0$ and all this process can be made in such a way that
 the nice cylinder is in $U_1^c.$
\B

\vspace*{-5mm}


%


\newpage
\begin{defi}(\textbf{Lateral border})
Given  a differentiable arc $\gamma$ and $r>0.$
The \emph{lateral border} $S$ of the cylinder $C(\gamma,r)$ 
is $\partial C(\gamma,r)$ minus 
 the top and bottom sides of the cylinder if they exist.
\end{defi}

Observe that nice cylinders have lateral borders.

%

\begin{defi}(\textbf{Separated horizontally})\label{def1.20} 
We say that a nice cylinder $C(\gamma,r)$ is \emph{separated
horizontally} by a set $\Lambda$ if there exists a connected
component of $\Lambda,$ let say $\Lambda_c,$ such that:
\begin{itemize}
\item $\Lambda_c$ intersects $C(\gamma,r)$ across the lateral border.
\item $C(\gamma,r)$ minus  $\Lambda_c$ 
has at least two connected component.
 \end{itemize}
\end{defi}

Now, we are going to prove that the locally maximal set for $f,$
found in section \ref{ss13}, has the geometrical
property of separating horizontally these nice cylinders as the one in claim \ref{afir6}.

%
%
%

\begin{lema}\label{afir3}
Given any arc $\gamma$ in $U_2^c$ with diameter greater than $\delta_0$
that admits a nice cylinder as in claim (\ref{afir6}), 
$\Lambda_f$ separates horizontally its nice cylinder.
\end{lema}

\pp
Let us denote by $T$ the nice cylinder associated to $\gamma$ as in the statement
and let
$A$ and $B$  denote the top and bottom sides of $T$ respectively. Let $\varepsilon>0$
be arbitrarily small.

\begin{figure}[h]
\begin{center}
\includegraphics[scale=0.35]{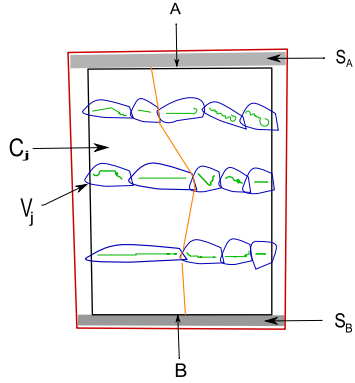}
\end{center}
\caption{\textit{$\Lambda_f$ splits "horizontally" every nice cylinder in at least two connected component}}\label{graf1.1b}
\end{figure}

Let $T'$ be a bigger cylinder containing $T$ joint together with two security regions,
denote by
$S_A$ and $S_B,$ and such that the distance between the lateral border of $T$
and the lateral border of $T'$ is small, for instance $d_H(T, T')=\frac{d_1}{6k}$, see figure (\ref{graf1.1b}).
For security regions $S_A$ and $S_B$ we mean two strips of $\frac{d_1}{6k}$ thickness
glued to the sides $A$ and $B$ of $T$,
or in other words, $S_A$ (respectively $S_B$) is the set of points in $T^c$ such that the distance from these points
to $A$ (respectively $B$) is less or equal to  $\frac{d_1}{6k}.$
This set $T'$ was constructed in such a way that its diameter is  greater
than $\delta'_0.$

Since $\gamma$ is in $U_2^c$ and its diameter is greater than $\delta_0,$ we can assure
that $T\cap \Lambda_f$ is non empty.  Consider all the connected components of $T\cap \Lambda_f.$
For every $x\in T\cap \Lambda_f,$ we assign  $\emph{K}_x$ the connected component
of $T\cap \Lambda_f$ that contains $x.$ Observe that we may define an equivalence relation:
$x\sim x'$ if and only if $\emph{K}_x=\emph{K}_{x'}.$
Then we pick one component from each class, or in other
words we pick just the connected components that are  two by two disjoints.

We claim that $\Lambda_f$ separates  $T$ horizontally, i.e.;
 there exists one component $\emph{K}_x$ such that $K_x\cap \partial T\neq \emptyset$ and $K_x$ separates $T$  in more than one connected
component.

Suppose it does not happen, i.e. none of the $\emph{K}_x$ separates
$T$ horizontally. Take $\emph{U}_x$ open set in $T'$ such that
$\emph{K}_x\subset \emph{U}_x,$ $\partial \emph{U}_x\cap \Lambda_f=\emptyset,$
$\partial \emph{U}_x$ is connected and $\partial \emph{U}_x$ does not divides  $T$ horizontally.
If there are many $\emph{K}_y$ accumulating in one $\emph{K}_x,$ then we could have a same
open set $\emph{U}_x$ containing more than one connected component $\emph{K}_y.$

Observe that the collection $\{ \emph{U}_x\}$ is an open cover of $T\cap \Lambda_f.$
Since it is compact, there is a finite subcover $\{ \emph{U}_i\}_{i=1}^N,$ i.e.
$T\cap \Lambda_f \subset \mathcal{U}:=\bigcup_{i=1}^N  \emph{U}_i.$

If the connected components of $\mathcal{U}$ does not separates horizontally
$T,$ it is easy to cons\-truct a curve going from $A$ to $B$ with
diameter greater than $\delta_0$ and empty intersection with the
$\emph{U}_i$'s; hence, this curve
does not intersects the set $\Lambda_f.$
But this contradicts the fact that every curve in $U_1^c$ with diameter
larger than $\delta_0$ intersects $\Lambda_f.$
 Then the connected components of $\mathcal{U}$  separate $T$ horizontally,
denote by $\emph{C}_j$ the connected components of $T$ minus these connected components of
$\mathcal{U}$ that separates $T$ horizontally.

Observe that every $\emph{C}_j$ is path connected, since
 they are the complement of a finite union of open sets in a simply connected set $T$.
 There exist a finite quantity of  $\emph{C}_j,$
 let us say $m$.
We can reorder these sets  enumerating from the top side.
If we denote by $V_j$ each of the connected components of $T\cap \mathcal{U}$ that separate
$T$ horizontally, we have two cases, either $\emph{C}_j$ is in between two consecutive
 $V_j$ and $V_{j+1}$ (or  $V_{j-1}$ and $V_j$)
or  $\emph{C}_j$ just intersects one $V_j$ on the border.

The idea is to build a curve from top to bottom of $T$  connecting
$\emph{C}_j$ with $\emph{C}_{j+1}$ in such a way that the diameter
of the arc is greater than $\delta_0$ but without intersecting
$\Lambda_f,$ which is a contradiction
because it is again in $U_1^c$ and  has diameter greater than
$\delta_0,$ then this curve must intersects $\Lambda_f.$

It is enough to show that we can pass from $\emph{C}_j$
to $\emph{C}_{j+1}$ without touching $\Lambda_f.$ For this, we must observe that
every $V_j$ is a union of finitely many  $\emph{U}_i$, let us say
$\emph{U}_{i_1},\ldots, \emph{U}_{i_j}.$
Pick a curve $\gamma_j$ in $\emph{C}_j$ going from top to bottom,
i.e.  $\gamma_j$ goes from $\partial V_j$ to $\partial V_{j+1}$
(or  $\partial V_{j-1}$ and $\partial V_j$)
and $\gamma_j$ does not intersects the interior of  $V_j$ and $V_{j+1}$
(or  $V_{j-1}$ and $V_j$), then
there exists $i_s\in \{i_1,\ldots,i_j\}$ such that
$\gamma_j\cap \partial\emph{U}_{i_s}\neq\emptyset.$ After that continue this arc
 picking a curve following
by the border of
$\emph{U}_{i_s}$ until $\emph{C}_{j+1},$ which has empty intersection with $\Lambda_f$
by construction,
if it is not possible to do in one
step, pick another $\emph{U}_{i_k}$ and repeat the process. Note that this
process finish in finitely many times.
The resulting arc from joint together all this segment
has diameter greater than $\delta_0$ and with
empty intersection with $\Lambda_f$ as we wanted.
\B

\begin{obs}\label{obs4}
In proposition \ref{afir2}, remembering that $d(h_g,id)<\eta,$  we may fix
$\eta<\min\{\frac{d_1}{6k},\delta_0,\beta\}.$
So for this $\eta,$
there exists $\varepsilon_0>0$ given by the shadowing lemma, see lemma \ref{sl} and proposition \ref{afir2},
and $\varepsilon_0$ determine  $\mathcal{V}_1(f)$ given in proposition \ref{afir2}.
\end{obs}

\begin{lema}\label{lema2}
Given $g\in \mathcal{V}_1(f)$ and given  an arc $\gamma$ in
$U_2^c$ with diameter greater than $\delta'_0$ such that it admits
 a  nice cylinder $C(\gamma,\frac{d_1}{2}),$
then $\gamma\cap \Lambda_g$ is not empty.
\end{lema}

\pp
Let $g\in \mathcal{V}_1(f).$ 
Take  an arc $\gamma$ in
$U_2^c$ with diameter greater than $\delta'_0$ such that
$C(\gamma,\frac{d_1}{2})$ is a  nice cylinder.

\begin{figure}[h]
\begin{center}
\includegraphics[scale=0.7]{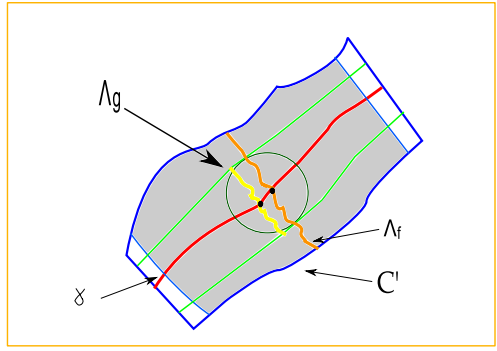}
\end{center}
\caption{\textit{$\Lambda_g$ intersects $\gamma$}}\label{graf1.4}
\end{figure}

By construction, we may assume that
every arc taken in the nice cylinder that goes from top to bottom has diameter
greater or equal to the diameter of $\gamma.$
We take  two security regions inside the cylinder,
 in the top and bottom sides of the cylinder respectively,
with $\frac{d_1}{6k}$ of thickness each one, i.e.
two strips
glued to the top and bottom sides of the cylinder such that
each one is  the set of points in the cylinder within distance
to top (respectively bottom) side less or equal to  $\frac{d_1}{6k},$ see figure (\ref{graf1.4}).
Let us denote by $C'$ the cylinder resulting of taking out
these two security strips from the original cylinder
$C(\gamma,\frac{d_1}{2}),$ then the diameter of $C'$ is still greater than $\delta_0$.

Hence, the diameter of $\gamma'=\gamma\cap C'$ is greater than $\delta_0$ and it is in $U_1^c.$
Lemma \ref{afir3} implies that $\Lambda_f$
separates horizontally $C',$ hence $\gamma'$ intersects $\Lambda_f,$ let us denote
by $x_f$ the point in the intersection.


 Since $x_f\in \Lambda_f,$ proposition \ref{afir2} and remark (\ref{obs4}), there exists
$x_g\in \Lambda_g\cap \mathbb{B}_{\eta}(x_f).$
Note that $\Lambda_f$ separates $\mathbb{B}_{\eta}(x_f)$ in at least two connected component.
Because $f\mid_U$ and  $g\mid_U$ are conjugate, follows that $\Lambda_g$ separates $\mathbb{B}_{\eta}(x_f)$ in at least two connected component as well.
Therefore, $\Lambda_g$ must intersects $\gamma.$\B

\vspace*{-5mm}

\subsection{Getting Arcs of Large Diameter}

In this subsection we show that under the hypothesis of volume expanding,
the diameter of the iterates of an open set grows on the covering.

\begin{lema}\label{afir4}
For every $g\in\mathcal{V}_1(f)$ and given $V$  an open path
connected set in $\mathbb{T}^n,$ there exists $m_0=m_0(V,g)\in
\mathbb{N}$ such that $diam(\tilde{g}^{m_0}(\widetilde{V}))>m,$
where $\tilde{g}$ and $\widetilde{V}$ are  the lift of $g$ and
$V,$ respectively, and $m$ was given in claim \ref{afir6}. In
particular, it contains an arc with diameter greater than $m$.
\end{lema}

\pp Let us suppose
that the Lemma is false in the covering space.
 Let $g\in\mathcal{V}_1(f)$
and $V$ be an open path connected set in $\mathbb{T}^n.$ If there
exists $k_0>0$ such that
$d_k=diam(\tilde{g}^{k}(\widetilde{V}))<k_0$, then
$vol(\tilde{g}^{k}(\widetilde{V}))\leq \left(
\frac{d_k}{2}\right)^n.$ But since $g$ is volume expanding, there
exists constant $\lambda>1$ such that
$vol(\tilde{g}^k(\widetilde{V}))>\lambda^k vol(\widetilde{V}),$
for $k\geq 1.$ 
Iterating by $\tilde{g}$, and since $\tilde g$ is a
diffeomorphisms in the covering space, the volume increase and
furthermore the diameter of its iterates grows  in the covering
space. Hence, there exists $m_0\in \mathbb{N}$ such that
$diam(\tilde{g}^{m_0}(\widetilde{V}))>m.$ \B

\begin{obs}\label{obs6}
For the case that $V$ is an open connected set, observe that
given a point in $V$ there exists an open ball centered in this point and contained
in $V$ such that it is path connected. Then we may apply Lemma \ref{afir4} to this ball and obtain a similar
statement for $V$.
\end{obs}

\subsection{Getting Sets of Large Radius}

In this subsection, we show that open sets intersecting $\Lambda_g,$
for $g$ close enough to $f,$ has large internal radius after large iterates.

\begin{lema}\label{lema3}
There exist $\mathcal{V}_2(f)$ and $R>0$ such that for every $g\in \mathcal{V}_2(f),$
if there is $x\in M$ such that
$g^n(x)\not\in U_0$ for every $n\geq 0,$ then there is $\varepsilon_0>0$
such that for every $0<\varepsilon<\varepsilon_0,$ there exists
$N=N(\varepsilon)\in\mathbb{N}$
such that $\mathbb{B}_R(g^N(x))\subset g^N(\mathbb{B}_{\varepsilon}(x)).$
\end{lema}

\pp
We may pick $U_3$ an open subset contained in $U_0$ such that
$m\{Df\mid_{U_3^c}\}>\lambda',$ with $1<\lambda'<\lambda_0.$
Take $\mathcal{V}_2(f)$ an open subset perhaps smaller than $\mathcal{V}_1(f)$ such that
 $m\{Dg\mid_{U_3^c}\}>\lambda'$ holds for every $g\in \mathcal{V}_2(f).$
Let us fix $R=d_H(U_0,U_3)>0.$

Given $0<\varepsilon<R,$ take $N\in\mathbb{N}$ such that
$(\lambda')^{-N}R<\varepsilon/2.$ Then
$\mathbb{B}_{(\lambda')^{-N}R}(x)\subset \mathbb{B}_{\varepsilon}(x).$

Observe that $\mathbb{B}_R(g^n(x))\cap U_3=\emptyset,$ for
every $n\geq 0.$ Also,
$$g^k(\mathbb{B}_{(\lambda')^{-N}R}(x))=\mathbb{B}_{(\lambda')^{-N+k}R}(g^k(x))
\subset \mathbb{B}_R(g^k(x)),$$ for every $0\leq k\leq N.$
In particular,
$g^k(\mathbb{B}_{(\lambda')^{-N}R}(x))\cap U_3=\emptyset,$ for
every $0\leq k\leq N.$
Then
$$g^N(\mathbb{B}_{(\lambda')^{-N}R}(x))=\mathbb{B}_R(g^N(x))\subset g^N(\mathbb{B}_{\varepsilon}(x)).$$
\B

\begin{obs}\label{obs5}
Let us note that  lemma \ref{lema3} holds for every point in $\Lambda_g.$
\end{obs}

%




\subsection{End of the Proof of  Main Theorem}\label{PMT}


Let $f\in E^1(\mathbb{T}^n)$ 
satisfying the hypotheses of the Main Theorem.
Lemma \ref{lema1} implies that we may assume the existence of
$\Lambda_f$ an expanding  locally maximal set for $f$.

Fix $0<\alpha<R,$ arbitrarily small. Given $x\in \mathbb{T}^n,$
since $\{ w\in f^{-i}(x): i\in\mathbb{N}\}$ is dense, there exists
$n_0\in\mathbb{N}$ such that
\begin{center}
$\{ w\in f^{-i}(x): 0\leq i \leq n_0\}$ is $\alpha/2$-dense.
\end{center}

Take a neighborhood $\mathcal{U}(f)\subset \mathcal{V}_2(f),$
where $\mathcal{V}_2(f)$ was given in
 lemma \ref{lema3}, such that for every $g\in \mathcal{U}(f)$
follows that
\begin{center}
$\{ w\in g^{-i}(x): 0\leq i \leq n_0\}$ is $\alpha/2$-close to $\{
w\in f^{-i}(x): 0\leq i \leq n_0\}$.
\end{center}

Hence, $\{ w\in g^{-i}(x): 0\leq i \leq n_0\}$ is $\alpha$-dense.

Let $V$ be an open connected set in $\mathbb{T}^n.$
By lemma \ref{afir4},
there exists $m_0\in \mathbb{N}$
such that $diam(\tilde{g}^{m_0}(\widetilde{V}))>m.$ Then we may pick an arc $\gamma$
in $\tilde{g}^{m_0}(\widetilde{V})$ with diameter larger than $m$ and applying claim (\ref{afir6})
follows that
there exists a connected piece $\gamma'$ of $\gamma$ such that
$\gamma^*=\pi(\gamma')$ is in  $U_2^c,$
diameter of $\gamma^*$ is larger than $\delta'_0$ and it admits a nice cylinder
$C(\gamma^*,\frac{d_1}{2}).$
By lemma \ref{lema2}
follows that $\gamma^*\cap\Lambda_g$ is not empty, let $y$ be a point in the intersection.

Hence, for this point $y,$ there exists $\varepsilon_0=\varepsilon_0(y)>0$
 such that $\mathbb{B}_{\varepsilon_0}(y)\subset g^{m_0}(V),$
by lemma \ref{lema3} taking $0<\varepsilon<\varepsilon_0,$ we get that
there exists $N=N(\varepsilon)\in\mathbb{N}$ such that
$$\mathbb{B}_R(g^N(y))\subset g^N(\mathbb{B}_{\varepsilon}(y))\subset g^{m_0+N}(V).$$
Hence, $\mathbb{B}_{\alpha}(g^N(y))\subset  g^{m_0+N}(V).$
Since the $\alpha-$density, we have that
$$\{ w\in g^{-i}(x): 0\leq i \leq n_0\}\cap \mathbb{B}_{\alpha}(g^N(y))\neq\emptyset.$$
Therefore, denoting by $p=m_0+N,$
$$\{ w\in g^{-i}(x): 0\leq i \leq n_0\}\cap g^p(V)\neq\emptyset.$$

Taking the $p-$th pre-image by $g$, we obtain that there is
$i_0\in\mathbb{N}$ such that
 $$\{ w\in g^{-i}(x): 0\leq i \leq i_0\}\cap V\neq\emptyset.$$

 Thus,   for every $g\in \mathcal{U}(f)$ follows that
 $\{ w\in g^{-i}(x): i\in\mathbb{N}\}$ is dense
 in $\mathbb{T}^n$ for every $x\in \mathbb{T}^n$.\B

\begin{figure}[h]
\begin{center}
\includegraphics[scale=0.4]{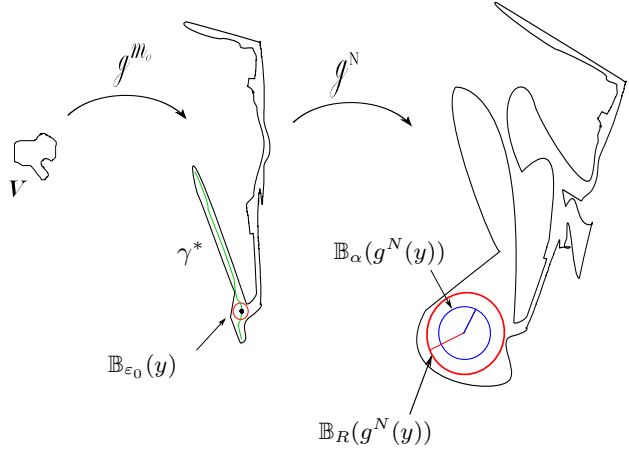}
\rput[0,0](-7,0.5){\small $\mathbb{B}_{\varepsilon_0}(y)$}
\rput[0,0](-4,2){\small $\mathbb{B}_{\alpha}(g^N(y))$}
\rput[0,0](-4.2,-0.4){\small $\mathbb{B}_R(g^N(y))$}
\rput[0,0](-6.2,2){ $\gamma^*$}
\end{center}
\caption{\textit{Iterations by the perturbed map}}\label{graf1.2}
\end{figure}

\subsection{The Main Theorem Revisited}\label{ss17}

In this section, we state a general  geome\-tri\-cal  version of the Main Theorem. 
Observe that using hypotheses (2) and (3) of the Main Theorem,
we showed in sections \ref{ss13} and \ref{ss15} the existence of a locally maximal
expanding set for $f$ which separates large nice cylinders, and in section \ref{ss14}
we proved that this geometrical property persist under perturbation, i.e.
there is a set $\Lambda_f$ locally maximal which intersects a nice class of arcs in $U_0^c$ 
and this pro\-per\-ty also holds for the perturbed. The hypothesis of
$f$ being volume expanding guarantees that given any open set in
the covering space, we are able to choose some iterates such that
it contains an arc with diameter large enough to apply claim
\ref{afir6} and lemma \ref{afir3}.  Hence, the Main Theorem may be
enunciated as follows: \vspace{6mm}

\begin{mtr}
\emph{Let  $f\in E^1(\mathbb{T}^n)$ be  volume expanding 
such that 
the pre-orbit of every point are dense. Suppose that there exist an open set $U_0$ with $diam(U_0)<1$ and
 $\Lambda_f$ a locally maximal expanding set for $f$ in $U_0^c$
such that  every  arc $\gamma$ in $U_0^c$ with
      diameter large enough  intersect $\Lambda_f$.
Then, the pre-orbit of every point are $C^1$ robustly dense.}
\end{mtr}
\vspace{2mm}

Observe that in the present version, it is already assumed the existences of an expanding
locally maximal invariant set that intersects large enough  arcs.
The proof goes showing  that the separation property is robust and this
is done in the same way that is done in the Main Theorem.

\begin{obs}
Note that the Main Theorem implies the Main Theorem Revisited,
but we do not know if  the reciprocal is true.
\end{obs}


\section{Robust transitive endomorphisms with invariant splitting} \label{ss121}
Now, we consider the case that the endomorphisms exhibits a type of partially hyperbolic splitting.
First we give the definition of partially hyperbolic endomorphisms which is slightly different than the
one for diffeomorphisms due to the fact that any point has different pre-images which implies that the
unstables subbundles are not unique (actually, they depend on the inverse branches).

\begin{defi}(\textbf{Unstable cone family})
Given $f:M\rightarrow M$ a local diffeomorphism, let $V$ be an open subset of $M$ such that
$f\!\!\mid_V$ is a diffeomorphism onto its image. Denote by $\varphi$  the inverse branches of $f$
restricted to $V;$ more precisely, $\varphi: f(V)\rightarrow V$ such that $f\circ\varphi (x)=x$
if $x\in f(V).$
A continuous cone field $\mathcal{C}^u=\{\mathcal{C}^u_x\}_{_x}$
defined on  $V$ is called \emph{unstable} if it is
forward invariant:
$$Df(x')\,\mathcal{C}^u_{x'}\subset \mathcal{C}^u_{f(x')}$$
for every $x'\in V\cap \varphi(V).$
\end{defi}

\begin{obs}
Given a point $x,$ there is not necessarily a unique unstable subbundle, i.e.
for each inverse path $\{x_k\}_{k\geq 0},$ it means $x_0=x$ and $f(x_{k+1})=x_k$ for
$k\geq 0,$ there exists an unstable direction belonging to $\mathcal{C}^u.$
\end{obs}

\begin{defi}(\textbf{Complementary splitting})
We say that a splitting $\mathbb{E}^c_x+\mathcal{C}^u_x$ is \emph{complementary}
if the unstable cone $\mathcal{C}^u_x$
 contains an invariant subspace whose dimension is
equal to the dimension of the manifold minus the dimension of the central subbundle.
\end{defi}

\begin{defi}(\textbf{Partially hyperbolic endomorphism with expanding 
                extremal direction})
It is said that  an endomorphism $f$ is
\emph{partially hyperbolic with expanding extremal direction}
provided 
for every $x\in M$
there exists a complementary splitting $\mathbb{E}^c_x+\mathcal{C}^u_x,$ where
$\{\mathcal{C}^u_x\}_{_x}$ is a  family of unstable cones,
and there exists $0<\lambda<1$ such that
for every inverse branches $\varphi$ of $f$ follows that
\begin{enumerate}
\item $\|D\varphi(x) \,v\|<\lambda,$ for all $v\in \mathcal{C}^u_x.$
\item $\|Df(x')\mid_{\mathbb{E}^{^c}(x')}\| \|D\varphi(x)v\|<\lambda,$ for all
        $v\in \mathcal{C}^u_{x},$ where $\varphi(x)=x', \;f(x')=x.$
\end{enumerate}
\end{defi}

\subsection{Theorem \ref{teo2}: Splitting Version}

 Now, we state a version of the Main Theo\-rem
for the case when the tangent bundle splits into two non-trivial  subbundles, one with an expanding
behavior and the other one with nonuniform behavior but  dominated by the expanding one.

\begin{teo}\label{teo2}
Let $f\in E^1(\mathbb{T}^n)$ be a locally diffeomorphism partially hyperbolic with expanding extremal direction
satisfying the following pro\-per\-ties:
\begin{enumerate}
\item $\{ w\in f^{-k}(x): k\in\mathbb{N}\}$ is dense for every $x\in \mathbb{T}^n.$
\item There exist $\delta_0>0,$  $\lambda_0>1$ and $k_0\in \mathbb{N}$
        such that for every  $x\in \mathbb{T}^n$, if  $\gamma$ is a disc 
        tangent to the unstable cone $\mathcal{C}_x^u$ with 
      internal diameter larger than  $\delta_0,$ there exists a point $y\in\gamma$ such that
      $m\{Df^i\mid_{\mathbb{E}^c({f^k(y)})}\}>\lambda_0^i,$ for all $i>0,$ for all $k>k_0.$
\end{enumerate}
Then, for every $g$ close enough to $f,$ $\{ w\in g^{-k}(x):
k\in\mathbb{N}\}$ is dense for every $x\in \mathbb{T}^n.$
\end{teo}

\subsection{Proof of Theorem \ref{teo2}}

The proof of Theorem \ref{teo2} is  similar
to the proof given in \cite{PS1},
where it is proved that any partially hyperbolic diffeomorphism satisfying
a  hypothesis like the one stated in Theorem \ref{teo2}
 and such that the
strong stable foliation is minimal, then the strong stable foliation is robustly minimal.
The key hypothesis in the statement of the main theorem in  \cite{PS1}  says that in any compact
piece of the unstable foliation, there exists a point such that the central
bundle  has  uniform expanding behavior  along the forward orbit, and this is exactly what we have.
The goal consists in proving  that this property is robust under perturbation.

Given a local diffeomorphism $f$ as in the statement of Theorem \ref{teo2},
we want to show that any small perturbation $g$ preserve the property of density
of the pre-orbit of any point. Our strategy is to prove that given any disc 
tangent
to the unstable cones for $g$ with  
large enough internal diameter   has a point such that
the central direction along the forward orbit by $g$ is uniformly expanding.

Observe that  given any open set, since we have
a direction that is indeed expanding, the diameter along the unstable direction
of the iterates growth. Then we are able to pick a disc
inside this iterate
such that the 
disc is tangent to the unstable cones with 
diameter large enough to apply
the last property.
Hence, there exists a point which its forward orbit is expanding in all direction,
then there is some iterate such that it contains a ball of a fix radius $\varepsilon.$

Since $g$ is close enough to $f,$ we have that the pre-orbit by $g$
are $\varepsilon-$dense. Therefore, given any open set, by the property of the unstable discs, 
there exists an iterate such that it intersects the pre-orbit by $g$ of any point. Thus, we conclude
the density of the pre-orbit of any point by the perturbation.

Moreover, the proof of Theorem \ref{teo2} can also be performed in the spirit of Main Theorem.
In fact, it is possible to show that
$$\bigcap_{l\geq 0} f^{-l}(\{x: m\{Df^n\mid_{\mathbb{E}^c({f^k(x)})}\}>\lambda_0^n,\;n>0,\;k>k_0\})$$
is an invariant expanding set such that separates unstable discs.
This provides a geo\-me\-tri\-cal interpretation.
\B

\subsection{Remarks About the Main Theorem and Theorem \ref{teo2}}


Observe that  in the Main Theorem we asked for large arcs to contain points such that its forward i\-te\-ra\-tions
remain in the expanding region.
The same is required in Theorem \ref{teo2} but just for large  unstable discs: there is a point there that
its forward iterates remain in ``an expanding region''for  the central bundle.

The main difference in their proof
arise from the fact  that in the version with splitting, since we know that we have uniform
expansion in one direction, any 
disc with 
internal diameter larger than $\delta_0$ and tangent to this direction,
growths up to  length $\delta_1>\delta_0$ in a bounded uniform time, independently of the disc. 
Note that we cannot guarantee that only assuming   volume expansion.

Observe that in the Main Theorem is not assumed
that $f$ does not have any splitting.
In fact, it could also be partially hyperbolic.
However, knowing in advance that the endomorphism is partially hyperbolic
then it is possible to get sufficient conditions for robust transitivity
weaker than the one required by the Main Theorem.

\section{$C^1$ Robust  transitivity and   volume expansion} \label{s4}


%

Before showing the relation between $C^1$ robust  transitivity and   volume expansion
(Theorem \ref{teo0}), let us recall some definitions that are involved in the statement

\begin{defi}\label{RT}
The set $\Lambda_f(U)=\bigcap_{n\in\mathbb{Z}}f^n(\overline{U})$
is $C^r$ \emph{robustly transitive} if
$\Lambda_g(U)=\bigcap_{n\in\mathbb{Z}}g^n(\overline{U})$ is
transitive for every endomorphism $g$ $C^r$ close enough to $f$,
where  $U$ is an  open set. It is said that a map $f$ is $C^r$
\emph{robustly transitive} if there exists a $C^r$ neighborhood
$\mathcal{U}(f)$ such that every $g\in \mathcal{U}(f)$ is
transitive.
\end{defi}

\begin{defi}\label{nosp}
We say that $f$ restricted to an invariant set $\Lambda$ has \emph{no dominated
splitting in a $C^r$ robust way} if there exists a  $C^r$ open neighborhood $\mathcal{U}(f)$
of $f$ such that for every $g\in\mathcal{U}(f)$ the tangent space $T\Lambda$
does not admit any dominated splitting.
\end{defi}


\begin{teo}\label{teo0}
Let $f\in E^1(M)$ be a local diffeomorphism and $U$ open set in $M$ such that
  $\Lambda_f(U)=\bigcap_{n\in\mathbb{Z}}f^n(\overline{U})$ is
$C^1$ robustly transitive set and it has no splitting in a $C^1$ robust way.
Then $f$ is volume expanding.
\end{teo}

\pp
The proof of this theorem is similar to the one of Theorem 4 in \cite[pp.361]{BDP},
nevertheless we include the main steps of the proof.

Let us consider $f\in E^1(M)$ a local diffeomorphism and denote by
$\Lambda_f(U)$ (nontrivial) $C^1$ robustly transitive  set for $f$
(note that $U$ could be the entire manifold). The idea of the
proof is to assume that $f$ is not volume expanding and show that
for every $\mathcal{U}(f)$ $C^1$ neighborhood of $f$ in $E^1(M),$
there exists $\psi\in\mathcal{U}(f)$ such that $\psi$ has
a sink and therefore $\psi$ cannot be transitive.

Suppose that $f$ is not volume expanding.  Since $f$ is onto, it cannot
be uniform volume contracting in the entire manifold, so there are
points in the manifold such that we have expansion, i.e. $1 \leq |det(Df^{k}(x))|$  for some $k\geq 0$,
but it does not expand too much, i.e. $|det(Df^{k}(x))|\leq 1+\epsilon,$ with $\epsilon$ small.
Then there are
sequences $x_n\in \Lambda_f(U),$ $k_n\in \mathbb{N}$ and $\tau_n>1,$
with $k_n\rightarrow\infty$ and $\tau_n\rightarrow 1^+,$ such that
$$1 \leq |det(Df^{k_n}(x_n))|<\tau_n^{k_n}.$$
This is equivalent to say that
$$\frac{1}{k_n}\sum_{i=0}^{k_n-1}\log(|det(Df(f^i(x_n)))|)<\log(\tau_n).$$
We may take $k_n$ such that $f^i(x_n)\neq f^j(x_n)$ for all
$i\neq j,$ $i,j\in\{0,\ldots,k_n\}.$ Consider for each $n$ the Dirac measure
$\delta_n$ supported in $\{x_n,f(x_n),\ldots,f^{k_n}(x_n)\},$ i.e.
$\delta_n=\frac{1}{k_n}\sum_{i=0}^{k_n-1}\delta_{f^i(x_n)}.$
As the space of probabilities is compact with the weak star topology,
there exists a subsequence of $\{\delta_n\}_n$ that converges
to an invariant probability measure $\mu$ such that
$$\int\log|det(Df(x))|d\mu(x)\leq 0.$$
In fact,  a classical argument proves that $\mu$ is invariant by $f,$ since
$f_*(\mu)-\mu$ is the weak star limit of $\frac{1}{k_{n_i}}(\delta_{f^{k_{n_i}}(x_{n_i})}-\delta_{x_{n_i}}),$
which converge to zero. Observing that
$$\begin{array}{ll}
  \int\log|det(Df(x))|d\delta_n &=\! \frac{1}{k_n}\!\sum_{i=0}^{k_n-1}\log(|det(Df^i(x_n))|)\!\\\\
&=\!\frac{1}{k_n}\log(|det(Df^{k_n}(x_n))|) \leq \log(\tau_n),
\end{array}$$
and since $\tau_n\rightarrow 1^+$ we deduce that  $$\int\log|det(Df(x))|d\mu(x) \leq 0.$$
By the ergodic decomposition theorem, there is  an ergodic and
$f-$invariant measure $\nu$ such that
$$\int\log|det(Df(x))|d\nu(x)\leq 0.$$

Using the ergodic closing lemma for nonsingular endomorphisms (see \cite{Castro}), given
$\varepsilon>0$ there is $g$ close to $f$ and a $g-$periodic point $y$
such that
$$\frac{1}{m_{\varepsilon}}\sum_{i=0}^{m_{\varepsilon}-1}\log(|det(Dg(g^i(y)))|)<\varepsilon,$$
where $m_{\varepsilon}$ is the period of $y.$ Note that if $\varepsilon\rightarrow 0,$ then
$m_{\varepsilon}\rightarrow\infty.$ So, taking $\varepsilon>0$ arbitrarily small and
$m_{\varepsilon}$ big,  using Franks' Lemma \cite{Franks} we get
$\varphi$ close to $g$ such that $\varphi^{m_{\varepsilon}}(y)=y\in\Lambda_{\varphi}(U)$ and
$$\frac{1}{m_{\varepsilon}}\sum_{i=0}^{m_{\varepsilon}-1}\log(|det(D\varphi(\varphi^i(y)))|)<0,$$
this means that
$|det(D\varphi^{m_{\varepsilon}}(y))|<\lambda<1.$ Observe that we
are assuming the dimension of the manifold
greater or equal to 2, so the fact that the modulus of the
jacobian of $\varphi$ be lower than 1 does not
imply that all the eigenvalues have modulus smaller than 1.

Since $\Lambda_{\varphi}(U)$ is $C^1$ robustly transitive, after
a perturbation, we may assume that the relative homoclinic class
$H(y,\varphi,U)$ of $y$ is the whole $\Lambda_{\varphi}(U)$ (see \cite{BDP} for definition).  
Now, consider the dense subset $\Sigma\subset\Lambda_{\varphi}(U)$
consisting of all the hyperbolic periodic points of $\Lambda_{\varphi}(U)$
homoclinically related to $y$.

If $\varphi$ is close enough to $f,$ then the tangent bundle does not admit
a splitting as well.  
Using the idea of the proof of Lemma 6.1 in \cite[pp. 407]{BDP}
and,  after that, Franks' Lemma,
we obtain that there exists $\psi$ a perturbation of
$\varphi$ and a point $p\in \Sigma$
such that all the eigenvalues of $D\psi^{m(p)}(p)$ have modulus strictly lower than
1, where $m(p)$ is the period of $p$. This means that the maximal invariant set
in $U$ of  $\psi$  contains a sink, but this is a contradiction since
 we choose $\psi$ sufficiently close to $f$
such that $\Lambda_{\psi}(U)$ is still transitive.
\B

\begin{obs}
If $\Lambda_f(U)$ admits a splitting, then the extremal indecomposable subbundle
is volume expanding. The proof is similar to the proof of theorem \ref{teo0}, restricting $Df$ to the extremal subbundle.
\end{obs}

\begin{obs}\label{obs16}
Theorem \ref{teo0} implies that  volume expanding of the extremal bundle
 is a necessary condition for an endomorphism, which is local diffeomorphism,
 to  be a robust transitive map. However, volume expanding is not a sufficient condition
 that guarantees robust transitivity
for a local diffeomorphism. For instance, consider a product of an expanding endomorphism
times an irrational rotation: this map
is volume expanding and transitive but not robust transitive.
\end{obs}

\begin{obs}
It is expected that if $f$ is robustly transitive and
has no invariant subbundles in a robust way,
then $f$ is a local diffeomorphism.
This result  depends on whether the Ergodic Closing Lemma holds
even if there are critical points, since for maps with critical
points already exists a version of Connecting Lemma, Closing Lemma and Franks' Lemma,
which are the principal results involved in the proof of Theorem \ref{teo0}.
\end{obs}

\section{Examples  of Robust Transitive Endomorphisms}\label{s5}

In this section we  show that there exist examples
of robust transitive endomorphisms verifying the hypotheses
of our main results. The first two examples co\-rres\-pond to endomorphisms without any
splitting where is applied the Main Theo\-rem and the revisited version, and they can 
be considered as an endomorphisms version of the one constructed in \cite{BV}; those examples
are \emph{Derived from Expanding} endomorphisms.  The last two ones correspond
to partially hyperbolic endomorphisms,  they can be considered  as an endomorphisms
version of the one constructed in \cite{BD} and \cite{NP}, and they are not isotopic to expanding ones.

\subsection{Example 1: Applying  Main Theorem}\label{ss31}

Consider $\mathcal{E}:\mathbb{T}^n\rightarrow\mathbb{T}^n$
an ex\-pan\-ding endomorphism, with $n\geq 2.$ Let us consider a Markov partition
of $\mathcal{E}$ and observe that its elements are given by $n-$dimensional closed rectangles.
Note that taking a large $m>1,$ the topological degree of
$\mathcal{E}^m$ is equal to the topological degree of $\mathcal{E}$
to the power $m-$th and the  Markov partition can be chosen in such a way that the number of its elements is equal to the topological degree of
 $\mathcal{E}^m,$ so without loss of generality  we may assume that the initial map
has many elements in the partition as we want. More precisely, if $N=topological\; degree(\mathcal{E})$,
we may assume that $N$ is large and therefore the Markov partition has $N$ elements.
Denote by $R_i$ the elements of the partition, with $1\leq i \leq N$;
 $R_i$ is closed, $int(R_i)$  is nonempty
and $int(R_i)\cap int(R_j)=\emptyset$ if $i\neq j$.

 Now, consider $\psi:\mathbb{T}^n\rightarrow\mathbb{T}^n$ a map isotopic to the identity
and denote by $\widehat{R}_i=\psi(R_i)$ for every $i.$ The idea of using this map
is to deform the elements of the Markov partition and get a new partition which elements are not
all of the same size (it could contain some very small elements and  others very big).

Set $U_0$ an open set in $\mathbb{T}^n$ such that if $\widetilde{U}$
 is the convex hull of the lift
of $U_0,$ then $\widetilde{U}\cap [0,1]^n$ is contained
in the interior of $[0,1]^n,$ i.e. $diam(U_0)<1.$
Note that there exists $\widehat{R}_i$ such that $\widehat{R}_i\cap U_0$ is nonempty. We also request that
 there are many $\widehat{R}_i$ contained in $U_0^c,$  observe that this condition is feasible
since the initial map has many elements in the partition.

\begin{figure}[h]
\begin{center}
\includegraphics[scale=0.5]{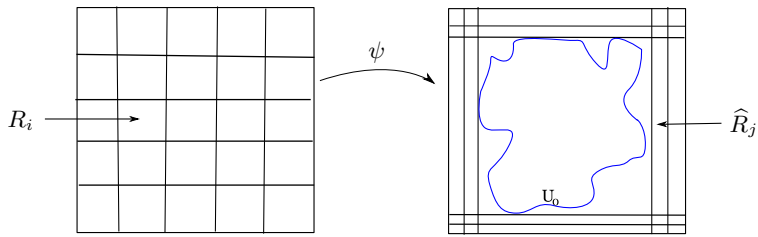}
\rput[0,0](-4.8,2.3){\small $\psi$}
\rput[0,0](-9.6,1.4){\small $R_i$}
\rput[0,0](0,1.3){\small $\widehat{R}_j$}
\end{center}
\caption{\textit{Deforming  the initial Markov Partition}}\label{graf3.1}
\end{figure}

Define $f_0:\mathbb{T}^n\rightarrow \mathbb{T}^n$ by
$f_0=\psi \circ \mathcal{E}.$ We assume that there
exist $p\in U_0$ and $q_i\in U_0^c$ fixed points for  $f_0,$ with $1\leq i\leq n-1.$ This is possible because
we may start with an expanding map which has as many fixed points as we need.

Let us suppose that $p$ and $q_i$ are expanding for $f_0$ in all directions, it means that
all the eigenvalues associated to these points has modulus greater than 1. Pick $\varepsilon>0$
small enough such that $\mathbb{B}_\varepsilon(q_i)\cap U_0=\emptyset$ and
$\mathbb{B}_\varepsilon(q_i)\cap \mathbb{B}_\varepsilon(q_j)=\emptyset$ if $i\neq j.$

Let us denote the decomposition of the tangent space as follows
$$T_x(\mathbb{T}^n)=\mathbb{E}_1^u\prec \mathbb{E}_2^u \prec \cdots\prec \mathbb{E}_{n-1}^u\prec \mathbb{E}_n^u,$$
where $\prec$ denotes that $\mathbb{E}_i^u$ dominates the expanding behavior of $\mathbb{E}_{i-1}^u.$

%


Next we deform $f_0$  by a smooth isotopy supported in $U_0\cup (\bigcup\mathbb{B}_\varepsilon(q_i))$  in such a way that:
\begin{enumerate}
\item The continuation of $p$ goes through a pitchfork bifurcation,
       appearing two periodic points $r_1, r_2,$ such that both are repeller
       and $p$ becomes a saddle point. But the new map $f$ still expand volume in $U_0.$
\item Two expanding eigenvalues of $q_i$ become complex expanding eigenvalues.
       More precisely, we mix the two expanding subbundles of $T_{q_i}(\mathbb{T}^n)$
       co\-rres\-pon\-ding to $\mathbb{E}_i^u(q_i)$ and $\mathbb{E}_{i+1}^u(q_i),$ obtaining
       $T_{q_i}(\mathbb{T}^n)=\mathbb{E}_1^u\prec \mathbb{E}_2^u \prec \cdots\prec
       \mathbb{F}_i^u\prec \mathbb{E}_{n-1}^u\prec\mathbb{E}_n^u,$ where $\mathbb{F}_i$ is two dimensional
       and correspond to the complex eigenvalues.
\item Outside $U_0\cup (\bigcup\mathbb{B}_\varepsilon(q_i)),$  $f$ coincides with $f_0.$
\item  $f$ is expanding in $U_0^c.$
\item There exists $\sigma>1$ such that $|det(Df(x))|>\sigma$ for every $x\in \mathbb{T}^n.$
\end{enumerate}

\begin{figure}[h]
\begin{center}
\includegraphics[scale=0.6]{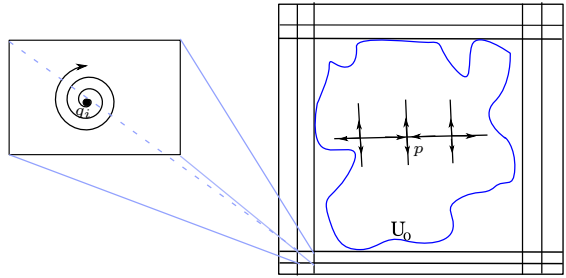}
\rput[0,0](-2.1,1.6){\tiny $p$}
\rput[0,0](-6.6,2.1){\tiny $q_i$}
\end{center}
\caption{\textit{$f$ isotopic to $f_0$}}\label{graf3.2}
\end{figure}

\subsubsection{Property of Large Arcs}

\begin{afir}
Every large arc in $U_0^c$ has a point such that its forward orbits remain in $U_0^c.$
\end{afir}

\pp
 Take $d$
the maximum of the diameter of the elements of the partition contained in $U_0^c.$
Note that every arc in $U_0^c$ with diameter larger
than $d$ cannot be contained in the interior of any element of the partition, more precisely
it has to intersect at least two elements of the partition. Hence, the image by $f$
of this arc $\gamma$ has diameter 1. So there exists a piece of $f(\gamma)$ in $U_0^c$
intersecting at least one element of the partition across two parallel sides, let us call
$\gamma^1$. Choose a pre-image of $\gamma^1$ in $\gamma$ and call it $\gamma_1.$

Repeating the process for $\gamma^1$, we have that there is $\gamma^2$ a piece of
$f(\gamma^1)$ verifying the same condition as $\gamma^1.$ Then, choose $\gamma_2$
a pre-image of $\gamma^2$ by $f^2$ in $\gamma.$

Thus, we construct a sequence of nested  arcs in $\gamma.$
The intersection is non empty and a point in this intersection satisfy our claim.
\B

\subsubsection{Remarks and variation of Example 1}

\begin{enumerate}
\item  $q_i$'s are fixed points for $f$  with complex expanding eigenvalues. Note that the e\-xis\-ten\-ce
        of these points ensures that the tangent bundle does not admit any invariant subbundle.
        We could also start with an expanding map having, besides $p,$ periodic points $q_i$
        with complex eigenvalues. In such a case, it is enough to make $p$ goes through a pitchfork bifurcation.
\item This example shows that $U_0$ can be as big as we desired while it verifies the hypothesis of having
     diameter less than 1.
\item It can be constructed in any dimension.
\end{enumerate}

\subsection{Example 2: Applying the Main Theorem Revisited} \label{ss32}

Let us consider $\mathcal{E}:\mathbb{T}^n\rightarrow\mathbb{T}^n$
an expanding endomorphism, with $n\geq 2.$ 
Assume that the initial map
has many elements in the  Markov partition, let us say $N$ elements.

\begin{figure}[h]
\begin{center}
\includegraphics[scale=0.2]{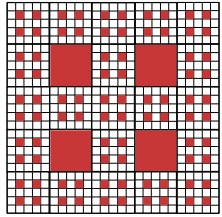}
\end{center}
\caption{\textit{$\{P_i\}_i$ collection}}\label{graf3.3}
\end{figure}

Denote by $R_i$ the elements of the partition, with $1\leq i \leq N$.
Since $\mathcal{E}$ is expanding, $R_i$ are closed, $int(R_i)$ are nonempty
and $int(R_i)\cap int(R_j)=\emptyset$ if $i\neq j$.
 Choose finitely many of these elements, $\{R_{i_j}\}_{j=1}^k,$ such that
 $R_{i_j}\cap R_{i_s}=\emptyset$ if $i_j\neq i_s,$ i.e. they are two by two disjoints.
 Consider the pre-images of every $R_{i_j},$ let us say
 $\mathcal{E}^{-1}(R_{i_j})=\{P_{i_j}^l\}_{l=1}^N.$
 Denote by $P_{i_j}^0=R_{i_j}.$ Next, we keep $P_{i_j}^r$ 
 such that
 $P_{i_j}^r\cap P_{i_{s}}^l=\emptyset$
 whenever $0\leq r\neq l\leq N$ and  $i_j\neq i_s.$
 Finally, let us denote by $\{P_i\}_i$ the collection of these latter subsets,
 so they are two by two disjoints. See figure \ref{graf3.3}.

\begin{figure}[h]
\begin{center}
\includegraphics[scale=0.3]{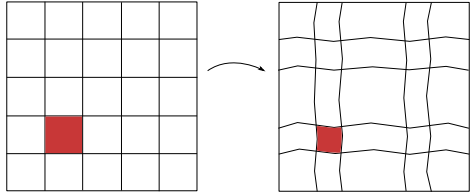}
\end{center}
\caption{\textit{Deforming the Markov partition}}\label{graf3.4}
\end{figure}

Now, consider $\psi:\mathbb{T}^n\rightarrow\mathbb{T}^n$ a map isotopic to the identity
and denote by $\widehat{P}_i=\psi(P_i)$ for every $i.$ See figure \ref{graf3.4}.

Choose $\widetilde{P}_i$ an open
connected subset such that its closure is contained in the interior of $\widehat{P}_i.$
Let $\phi_i:\mathbb{T}^n\rightarrow\mathbb{T}^n$ be a map isotopic to the identity
such that
\begin{itemize}
\item $\phi_i\mid_{\widetilde{P}_i}$ is not expanding.
\item $\phi_i\mid_{\widehat{P}_i^c}$ is the identity.
\end{itemize}

\vspace*{-5mm}
\begin{figure}[h]
\begin{center}
\includegraphics[scale=0.3]{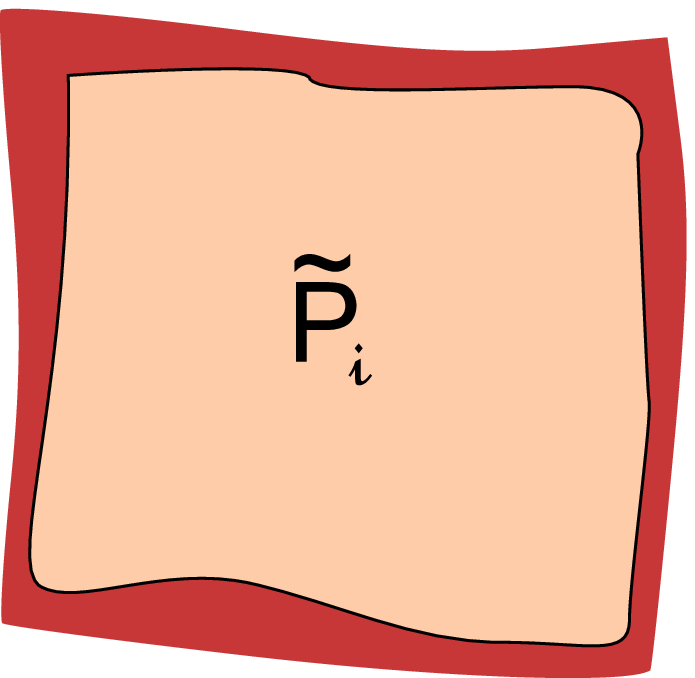}
\end{center}
\end{figure}

Define $\phi:\mathbb{T}^n\rightarrow\mathbb{T}^n$ by
$$\phi(x)=\left\{\begin{array}{ll}
                \phi_i(x),\;& \mbox{if}\; x\in\widehat{P}_i\\\\
                x,\;& \mbox{if}\; x\not\in \bigcup_i\widehat{P}_i
                \end{array}
                \right.$$
Hence, $\phi$ is equal to the identity in  $[\bigcup_i\widehat{P}_i]^c,$
 expands volume but is not expanding
in $\bigcup_i\widehat{P}_i.$

Once we have defined all these maps, we consider the map
$f=\phi \circ\psi \circ \mathcal{E}$ from $\mathbb{T}^n$ onto itself and denote
by $U_0=int (\bigcup_i\widehat{P}_i).$ Observe that $f$ verifies that:
\begin{enumerate}
\item[(\emph{i})] $f$ is a volume expanding endomorphism.
\item[(\emph{ii})] $f$ is an expanding map in $U_0^c.$
\item[(\emph{iii})] $\Lambda_f=\bigcap_{n\geq 0}f^{-n}(U_0^c)$ is an expanding locally maximal set for $f$
        which has the property that separate large nice cylinders.
\end{enumerate}

%


Since (\emph{i}) and (\emph{ii}) are immediate from the construction of $f,$ we concentrate
our interest in proving (\emph{iii}).

\subsubsection{$\Lambda_f$ Separates Large Nice Cylinders}

\begin{figure}[h]
\begin{center}
\includegraphics[scale=0.25]{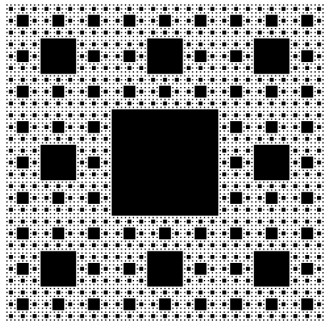}
\end{center}
\caption{\textit{$\Lambda_f$ looks like a Sierpinski set}}\label{graf3.3b}
\end{figure}

Note that by the construction of $U_0,$ we have that the elements of the  pre-orbit
of $U_0$ are two by two disjoints.
Let us consider $d_0=\max\{diam(c.c. \bigcup_{k\geq 0} f^{-k}(U_0))\}.$
Since the definition of $U_0,$ $0<d_0<1.$ Observe that $\Lambda_f$ looks like  a Sierpinski set, see figure \ref{graf3.3b}.

\begin{afir}\label{afir31}
If $\gamma$ is an arc in $U_0^c$ with diameter $1,$ then
$\gamma$ intersects $\Lambda_f.$
\end{afir}

\pp
Let $\gamma$ be an arc in $U_0^c$ such that $diam(\gamma)=1.$
Suppose that $\gamma$ does not intersect $\Lambda_f.$

Remember that $\Lambda_f=\mathbb{T}^n\setminus \bigcup_{k\geq 0} f^{-k}(U_0),$
it means that if $x\in \Lambda_f,$ then $f^k(x)\not\in U_0$ for all $k\geq 0.$ Therefore, 
$\gamma$ is contained in one pre-image of $U_0$ or in a union of pre-images of $U_0.$

Observe that $\gamma$
cannot be contained in just one pre-image of $U_0,$ because if it is contained
in $f^{-k}(U_0)$ for some $k\geq 0$,
then $diam(\gamma)<diam(f^{-k}(U_0))<d_0,$ which is absurd because $d_0<1.$

Hence, $\gamma$ should be contained in a union of pre-images of $U_0,$ since $\gamma$ is
compact we can cover with a finite union of pre-images of $U_0$. But we know that
the pre-images of $U_0$ are two by two disjoints, hence there exist points in $\gamma$ that
cannot be covered by the pre-images of $U_0.$ In particular, $\gamma$ intersects $\Lambda_f.$
\B

\begin{obs}
We have already proved
the existence of the invariant expanding locally maximal set $\Lambda_f.$ Moreover,
by claim (\ref{afir31}) we get that this invariant set intersects every arc with large diameter.
Then by the Main Theorem Revisited follows that this map is robustly transitive.
\end{obs}

\subsubsection{Remarks About Example 2}
\begin{enumerate}
\item We can apply our Main Theorem Revisited to this example, obtaining in particular
        that $f$ is robustly transitive.
\item The $\widehat{P}_i$'s can be as many and  as big as we want.
\item We can construct many examples starting with this initial map.
         In particular, we can construct examples
     without invariant subbundles, such as putting a fix point in the complement
     of the $U_0$ with complex eigenvalues and doing
     a derived from an expanding endomorphisms inside of some $\widehat{P}_i$.
\end{enumerate}

\subsection{Example 3: Applying  Theorem \ref{teo2}} \label{ss33} 

The idea of next example is to build an endomorphism in the 2-Torus which
is a skew-product and contains  a ``blender" for endomorphisms. This example is more
or less a standard adaptation for endomorphisms of examples obtained in \cite{BD}
for diffeomorphisms. The main goal is to get ``blenders"  for endomorphisms and since
we do not necessarily need to deal with stable foliation, the task is easier than the case of
diffeomorphisms. For more information about blenders see \cite{BD}.

First,  let us identify the 2-Torus with $[0,1]^2$ and let us
establish some notation before defining the map. 
Pick $0\!\!<\!\!a\!<\!\!b\!<\!\!3/4$ and $1/4\!\!<\!\!c\!<\!\!d\!<\!\!1.$
Denote $J_1\!\!=\!\![0,b],\; J_2\!\!=\!\![a,3/4],
\;J_3\!\!=\!\![1/4,d]$ and $J_4\!\!=\!\![c,1].$ Note that $J_1\cap J_2=[a,b]$ and
$J_3\cap J_4=[c,d].$
 This decomposition is associated to
the horizontal fibers.

Next, fix $N>3$ and pick  $0<a_1<b_1<a_2<b_2<a_3<b_3<a_4<b_4<1$
such that  $b_i-a_i=1/N.$ Let us denote
 by $I_i=[a_i,b_i]$ with $1\leq i\leq 4.$
Note that  they are two by two disjoint
and do not contain 0 or 1.  We associate this decomposition
 to the vertical fibers.

Let us call $R_i=J_i\times I_i$ where $i=1,2,3,4.$

Now, define $\Phi:\mathbb{T}^2\rightarrow \mathbb{T}^2$ by
$$\Phi(x,y)=(\varphi_y(x), \mathcal{E}(y)),$$ where
$\varphi_y, \mathcal{E}: S^1\rightarrow S^1$ are defined
as follows:
\begin{enumerate}
\item $\mathcal{E}$ is an expanding endomorphism such that:
\begin{itemize}
  \item $\mathcal{E}(I_i)=[0,1]$ for every $i.$
  \item There exists $a_i<c_i<b_i$ such that  $\mathcal{E}(c_i)=c_i.$
\end{itemize}
\item $\varphi_y$ is defined by
      $\varphi_y(x)=f_i(x),$ if  $y\in I_i,$ where
      $f_i: S^1\rightarrow S^1$ are diffeomorphisms 
      defined as follows:
      \begin{itemize}
        \item $f_1$ and $f_2$ satisfy the following properties:
                \begin{enumerate}
                \item[(\emph{i})] $f_1$ has two fixed points, $0$ and $a'\!\!\in\!\!(3/4,1),$ where
                                    $0$ is a repeller and  $a'$ is an attractor for $f_1.$
                \item[(\emph{ii})] $f_2$ has two fixed points, $3/4$ and $a''\!\!\in\!\! (a',1),$ where 
                                    $3/4$ is a repeller and  $a''$ is an attractor  for $f_2.$
                \item[(\emph{iii})]$f_1(J_1)= f_2(J_2)=[0,3/4].$
                \item[(\emph{iv})] $|f'_i \mid_{J_i}|>1$ for $i=1,2.$
              \end{enumerate}
         \item $f_3$ and $f_4$ satisfy the following properties:
                \begin{enumerate}
                \item[(\emph{i}')] $f_3$ has two fixed points, $c'\!\!\in\!\!(0,1/4)$ and $1/4,$ where 
                                    $c'$ is an attractor  and  $1/4$ is a repeller for $f_3.$
                \item[(\emph{ii}')] $f_4$ has two fixed points, $c''\!\!\in\!\!(0,c')$ and
                                    $1,$ where
                                    $c''$ is an attractor and  $1$ is a repeller for $f_4.$
                \item[(\emph{iii}')] $f_3(J_3)= f_4(J_4)=[1/4,1].$
                \item[(\emph{iv}')] $|f'_i\mid_{J_i}|>1$ for $i=3,4.$
                \end{enumerate}
      \end{itemize}
\item $|det(D\Phi)|=|\frac{\partial \varphi_y}{\partial x}\,\mathcal{E}'|>1.$
\item $\mathcal{E}'\gg \frac{\partial \varphi_y}{\partial y}.$
\end{enumerate}

\begin{figure}[h]
\begin{center}
\includegraphics[scale=0.65]{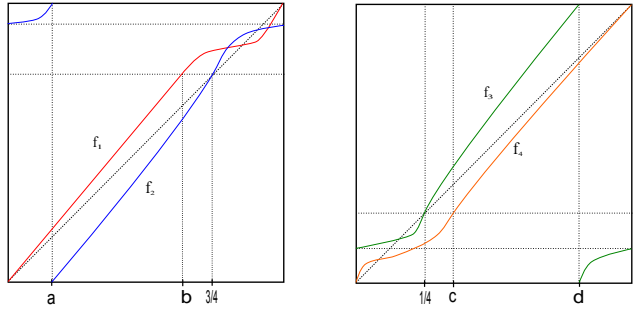}
\end{center}
\caption{\textit{Horizontal dynamics}}\label{ej3_1}
\end{figure}

\begin{figure}[h]
\begin{center}
\includegraphics[scale=0.75]{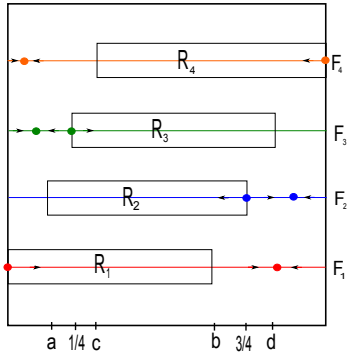}
\end{center}
\caption{\textit{This is how the dynamics $\Phi$ looks like}}\label{ej3_2}
\end{figure}

Hence, the horizontal fibers $F_i=S^1\times \{c_i\}$ are invariant
by  $\Phi,$ see figure \ref{ej3_2}. Moreover, by condition (4),
the image by $\Phi$ of every vertical fiber is almost a vertical fiber, in
the sense that the tangent vector is close to a vertical one; more
precisely, the unstable cones family are almost vertical.

Next, we consider $\Lambda_1^+=\bigcap_{n\geq 0}\Phi^{-n}(R_1\cup
R_2)$ and $\Lambda_2^+=\bigcap_{n\geq 0}\Phi^{-n}(R_3\cup R_4).$
Let $\Lambda_1=\bigcap_{n\in \mathbb{Z}}\Phi^{-n}(R_1\cup R_2)$
and $\Lambda_2=\bigcap_{n\in \mathbb{Z}}\Phi^{-n}(R_3\cup R_4).$
Note that both sets, $\Lambda_1$ and $\Lambda_2$ are expanding 
locally maximal invariant sets 
and each one contains a blender.

\subsubsection{$\Lambda_1$ and $\Lambda_2$ Separate Large Vertical Segments}

Let us denote by $\ell_1^u(p)$ the vertical segment passing through $p$ and length 1.

\begin{afir}\label{afir33}
For every $p\in R_1\cup R_2,$ follows that $\ell_1^u(p)\cap \Lambda_1^+\neq\emptyset.$
\end{afir}

\pp
Let $p\in R_1\cup R_2,$ then $L_i=\ell_1^u(p)\cap R_i$ is non empty for some $i\in\{1,2\}.$ 
The image of $L_i=\ell_1^u(p)\cap R_i$
by $\Phi$ has length 1 and by property (4) of $\Phi$ follows that $\Phi(L_i)$ is almost vertical.
Moreover, $L_i\cap F_i\neq\emptyset$ and $\Phi(L_i\cap F_i)\in F_i\subset R_i.$ 
Then $\Phi(\ell_1^u(p))\cap (R_1\cup R_2)\neq\emptyset.$
Call $K_1^i$ the connected component $\Phi(\ell_1^u(p))\cap R_i.$ 
Note that $P_2(K_1^i)=I_i,$ where $P_2$ is the projection in the second coordinate.
Consider the pre-image of $K_1^i$ by $\Phi$ in $L_i$
and call it $S_1^i.$

Now, iterate $K_1^i$ by $\Phi$, doing a similar process we obtain $K_2^i,$ the
connected component  $\Phi(K_1^i)\cap R_i$
such that $P_2(K_2^i)=I_i.$ Again take a pre-image of $K_2^i$ by $\Phi^2,$ giving a compact
segment $S_2^i\subset S_1^i.$
Repeating this process, we may construct a nested sequence of
compact segment $\{S_k^i\}_k$ in each $R_i.$ Thus, $\bigcap_k S_k^i$ is not empty
and belong to $\ell_1^u(p)\cap \Lambda_1^+.$\B

\begin{afir}\label{afir35}
For every $p\in R_1\cup R_2,$ follows that $\ell_1^u(p)\cap \Lambda_1\neq\emptyset.$
\end{afir}

\pp
By claim (\ref{afir33}), we know that there exists a point $z\in \ell_1^u(p)\cap \Lambda_1^+,$
this means that $\Phi^n(z)\in R_1\cup R_2$ for every $n\geq 0.$

Then, just remain to show that there exists a sequence
$\{z_k\}_{k\geq 0}\subset R_1\cup R_2$ such
that $z_0=z$ and $\Phi(z_k)=z_{k-1}.$ The idea of the construction
of such a sequence
 is to use now the property (2-\emph{iii}) of overlapping
in the horizontal dynamics.

Knowing that $\Phi(R_1)=f_1(J_1)\times [0,1]$
and $\Phi(R_2)=f_2(J_2)\times [0,1],$ since property (2-\emph{iii}) we get that
$\Phi(R_1)=\Phi(R_2)= [0,3/4]\times [0,1].$
Hence,  $z_0\in (R_1\cup R_2)\cap \Phi(R_1)$ or $z_0\in (R_1\cup R_2)\cap \Phi(R_2),$
then there exists $z_1\in R_1\cup R_2$ such that $\Phi(z_1)=z_0.$ Repeating this process
we construct the require sequence.
\B

\begin{afir}\label{afir34}
For every $p\in R_3\cup R_4,$ follows that $\ell_1^u(p)\cap \Lambda_2\neq\emptyset.$
\end{afir}

\pp
The proof is similar to claim (\ref{afir35}) just making the necessary arrangement.\B

\begin{afir}
For every $q\in\mathbb{T}^2,$ we have that either $\ell_1^u(q)\cap \Lambda_1\neq\emptyset$
or $\ell_1^u(q)\cap \Lambda_2\neq\emptyset.$
\end{afir}

\pp Given any point $q\in\mathbb{T}^2,$ note that $\ell_1^u(q)\cap
R_i\neq\emptyset$ for some $1\leq i\leq 4.$
Hence, taking $p_i\in\ell_1^u(q)\cap R_i$ and noting that
$\ell_1^u(p_i)=\ell_1^u(q)$, we may use claim (\ref{afir35}) or
(\ref{afir34}) to conclude that either $\ell_1^u(q)\cap
\Lambda_1\neq\emptyset$ or $\ell_1^u(q)\cap
\Lambda_2\neq\emptyset.$\B

\subsubsection{Remarks About Example 3}

This example was constructed in the 2-Torus with one dimensional central bundle,
 but we can construct it in any $\mathbb{T}^n$
and the dimension of the central bundle not need to be 1.
Also, we can use more than 4 dynamics in
the horizontal, that is more than four maps in the first variable.
More precisely, we put 2 blenders in the dynamic, induced by these
four maps, but we can consider as many blenders as we want.

\subsection{Example 4: Applying Theorem \ref{teo2}}\label{ss34}

Let $\mathbb{B}_0$ be an open ball in $\mathbb{T}^m$ centered at $0$ with radius $\alpha<1$  and
$\varphi_0:\mathbb{T}^m\rightarrow \mathbb{T}^m$ be a differentiable map isotopic to the identity
such that:
\begin{itemize}
\item $\varphi_0(0)=0$
\item There exist $0<\lambda_0<\lambda_1<1$ such that
        $\lambda_0<m\{D\varphi_0\}<|D\varphi_0\mid_{\mathbb{B}_0}|<\lambda_1,$ i.e.
        $\varphi_0$ is a contraction in a disk.
\end{itemize}

Let us consider $\mathbb{D}_0$ the lift of $\mathbb{B}_0$
to $\mathbb{R}^m$ and $\widetilde{\varphi}_0$
the lift of $\varphi_0.$ Note that $\widetilde{\varphi}_0(0)=0$ and
$\lambda_0<m\{D\widetilde{\varphi}_0\}<|D\widetilde{\varphi}_0\mid_{\mathbb{D}_0}|<\lambda_1.$
By Proposition 2.3  \cite{NP},
there exists $k\in\mathbb{N}$ such that for every small $\varepsilon>0,$
there exist $c_1,\ldots,c_k\in\mathbb{B}_{\varepsilon}(0)$ and $\delta>0$
such that $\mathbb{B}_{\delta}(0)\subset \overline{Orbit_{\mathcal{G}}^+(0)},$
where $\mathcal{G}=\mathcal{G}(\widetilde{\varphi}_0, \widetilde{\varphi}_0+c_1,\ldots, \widetilde{\varphi}_0+c_k)$
and $Orbit_{\mathcal{G}}^+(0)$ is the set of points lying on some orbit of $0$
under the iterated function system (IFS) $\mathcal{G};$ more precisely,
if we denote by $\widetilde{\phi}_0=\widetilde{\varphi}_0$ and $\widetilde{\phi}_i=\widetilde{\varphi}_0+c_i$ for $i=1,\ldots,k,$
 then $Orbit_{\mathcal{G}}^+(0)$ is the set
of sequence $\{\widetilde{\phi}_{\Sigma_l}(0)\}_{l=1}^{\infty}$ where  $\Sigma_l=(\sigma_1,\ldots,\sigma_l),$
$\widetilde{\phi}_{\Sigma_l}=\widetilde{\phi}_{\sigma_l}\circ\cdots\circ\widetilde{\phi}_{\sigma_1}$
 and $\{\sigma_i\}_{i\in\mathbb{N}} \in \{0,\ldots,k\}^{\mathbb{N}}.$
 (For more details about IFS see \cite{NP})

Now choose
$p_1,\ldots,p_r\in \mathbb{T}^m$ such that
$\mathbb{T}^m\subset \bigcup_j \mathbb{B}_{\delta}(p_j).$

If $\phi_i$ is the projection of $\widetilde{\phi}_i$ on $\mathbb{T}^m,$
define for every $j$ the IFS $\mathcal{G}_j=\mathcal{G}_j(\phi_0+p_j,\phi_1+p_j,\ldots,\phi_k+p_j).$
Then $\mathbb{B}_{\delta}(p_j)\subset \overline{Orbit_{\mathcal{G}_j}^+(0)}.$
Therefore, there exists an open set $\mathbb{D}_0\subset \mathbb{B}_0$
such that $\bigcup \phi_i(\mathbb{D}_0)\supset \mathbb{D}_0,$ i.e. the IFS has the covering property.
Hence, $\bigcup_i \phi_i(\mathbb{B}_{\delta'} (p_j))\supset \mathbb{B}_{\delta'} (p_j),$
with $0<\delta'\leq \delta.$ Moreover, $\mathcal{G}_j$ has also the overlapping property as in Example 3,
in the previous section.

Define the skew-product $F:\mathbb{T}^m\times \mathbb{T}^n \rightarrow \mathbb{T}^m\times \mathbb{T}^n$ by
$$F(x,y)=(\psi_y(x),\mathcal{E}(y)),$$
where:
\begin{itemize}
\item $\mathcal{E}:\mathbb{T}^n\rightarrow \mathbb{T}^n$ is an expanding map with $(k+1)r$ fixed points,
        let us denote the fixed points by $e_1^i,\ldots,e_r^i$ with $0\leq i\leq k.$
\item For every $y\in\mathbb{T}^n,$ $\psi_y:\mathbb{T}^m\rightarrow\mathbb{T}^m$
        is a differentiable map isotopic to the identity such that
        $\psi_{e_j^i}=\phi_i+p_j,$ with $0\leq i\leq k$ and $1\leq j\leq r.$
\end{itemize}

Hence, every fiber $\mathbb{T}^m\times\{e_j^i\}$ is invariant by $F.$
Set $R_j^i=\mathbb{B}_{\delta'} (p_j)\times Q_j^i,$ where  $Q_j^i$ is a small neighborhood
of $e_j^i$ in $\mathbb{T}^n$ such that $\mathcal{E}(Q_j^i)=\mathbb{T}^n$ and they are all disjoints
for every $i,j.$ Note that $R_j^i$ are the analogous of $R_i$ in the previous example.

Let $\Lambda_F:= \displaystyle\bigcap_{n\in \mathbb{Z}}F^n(\bigcup_{i,j}R_j^i).$

\subsubsection{$\Lambda_F$  Separate Large Unstable Discs}

\begin{afir}
$\Lambda_F$ verifies that for every $z\in\bigcup_{i,j}R_j^i$ follows that
$\ell_1^u(z)\cap \Lambda_F\neq\emptyset,$ where $\ell_1^u(z)$ is an
unstable disc of internal diameter 1 passing through $z.$
\end{afir}

\pp
 We may prove
that there exists a point $z\in \bigcup_{i,j}R_j^i$ such that $F^n(z)\in \bigcup_{i,j}R_j^i$ for
every $n\geq 0$ in a similar way as we proved claim (\ref{afir33}) in previous example.

Moreover, for this $z$ there exists $z_1\in \bigcup_{i,j}R_j^i$ such that $F(z_1)=z.$ In fact,
the idea is more or less the same as in previous example, we must note that
$F(R_j^i)=\psi_{e_j^i}(\mathbb{B}_{\delta'} (p_j))\times \mathcal{E}(Q_j^i)
=\phi_i(\mathbb{B}_{\delta'} (p_j))\times \mathbb{T}^n.$

On the other hand, using the property of covering and overlapping
follows that
$$\bigcup_{i,j}R_j^i= \bigcup_{i,j} \mathbb{B}_{\delta'} (p_j)\times Q_j^i\subset
\bigcup_{i,j} \phi_i(\mathbb{B}_{\delta'} (p_j)) \times \mathcal{E}(Q_j^i)= F(\bigcup_{i,j}R_j^i).$$
Therefore, since $z\in \bigcup_{i,j}R_j^i$, there exists $z_1\in \bigcup_{i,j}R_j^i$ such that
$F(z_1)=z.$ Inductively we can construct a sequence
$\{z_k\}_{k\geq 0}\subset \bigcup_{i,j}R_j^i$ such that $z_0=z$
and $F(z_k)=z_{k-1}.$

Thus, $z\in \ell_1^u(z)\cap \Lambda_F.$\B

\subsubsection{Remarks About Example 4}

This example is a generalization of Example 3. The intention here is to show that we may apply
 Theorem \ref{teo2} without taking into account the dimension of the central bundle and this could be as large as we want.
Another observation is that the existence of blenders guarantee that our examples are robustly transitive
and this example satisfies the property over the unstable discs with sufficiently large internal diameter
 intersecting the invariant
expanding locally maximal set for the skew-product.

\nocite{PS2}
\nocite{Katok}

\bibliographystyle{alpha}
\bibliography{library1}

%
%
%
%

\end{document}